\documentclass[journal ]{new-aiaa}
\usepackage[utf8]{inputenc}
\usepackage{textcomp}

\usepackage{fullpage,graphicx,psfrag,amsmath,amsfonts,verbatim,mathtools,dsfont,xcolor,hhline,algorithmicx,algorithm,algpseudocode,array,tikz,hyperref,leftindex,listings,bm,subcaption,pgfplots,booktabs,calc}
\pgfplotsset{compat=1.18}
\usepgfplotslibrary{fillbetween}

\usepackage[version=4]{mhchem}
\usepackage{siunitx}
\usepackage{longtable,tabularx}
\newcommand{\ones}{\mathbf 1}

\newcommand{\diag}{\mathop{\bf diag}}

\newcommand{\norm}[1]{\left\lVert#1\right\rVert}
\newcommand{\mnorm}[1]{{\left\vert\kern-0.25ex\left\vert\kern-0.25ex\left\vert #1 
    \right\vert\kern-0.25ex\right\vert\kern-0.25ex\right\vert}}

\newcommand{\argmin}{\mathop{\rm argmin}}

\newcommand{\blkdiag}{\mathop{\rm blkdiag}}

\newcommand{\eg}{{\it e.g.}}
\newcommand{\ie}{{\it i.e.}}

\definecolor{proxrandom_color}{rgb}{0,0,1}
\definecolor{proxws_color}     {rgb}{1,0,0}
\definecolor{ipoptran_MS_color}{rgb}{0.4940,0.1840,0.5560}
\definecolor{ipoptws_MS_color} {rgb}{0.3010,0.7450,0.9330}
\definecolor{ipoptran_SS_color}{rgb}{0.9290,0.6940,0.1250}
\definecolor{ipoptws_SS_color} {rgb}{0,1,0}

\definecolor{color1}{rgb}{1.0000,0.0000,0.0000}   
\definecolor{color2}{rgb}{0.0000,0.0000,1.0000}   
\definecolor{color3}{rgb}{0.0000, 0.7000, 0.0000}   
\definecolor{color4}{rgb}{0.0000, 0.6000, 0.6000}   
\definecolor{color5}{rgb}{0.7000, 0.3000, 0.7000}  
\definecolor{color6}{rgb}{0.5500,0.0500.0.0500}   
\definecolor{color7}{rgb}{0.8500,0.3250,0.0980}
\definecolor{color8}{rgb}{0,0,0}

\pgfplotsset{every axis plot}
\pgfplotsset{grid style=dotted}

\title{Sequential Convex Programming with Filtering-Based Warm-Starting for Continuous-Time Multiagent Quadrotor Trajectory Optimization}

\author{Minsen Yuan\footnote{Graduate Research Assistant, Department of Aerospace Engineering and Mechanics. Corresponding author. Email: \texttt{yuan0450@umn.edu}.}
       and Yue Yu\footnote{Assistant Professor, Department of Aerospace Engineering and Mechanics, AIAA Member. Email: \texttt{yuey@umn.edu}.}}
\affil{University of Minnesota Twin Cities, Minneapolis, MN 55455.}

\begin{document}

\begin{tikzpicture}
  \begin{axis}[
    hide axis,                 
    xmin=0, xmax=1,            
    ymin=0, ymax=1,            
    width= \linewidth,
    height = 0.3\linewidth, 
    legend to name=sharedlegend_2,
    legend style={
      legend columns=6,
      /tikz/column 2/.style={column sep=20pt},
      /tikz/column 4/.style={column sep=20pt},
      /tikz/column 6/.style={column sep=20pt},
      /tikz/column 8/.style={column sep=20pt},
      /tikz/column 10/.style={column sep=20pt},
      font=\footnotesize
    },
    legend image post style={line width=1.5pt}
  ]
    \addlegendimage{line width=1pt, color=color1}\addlegendentry{Agent 1}
    \addlegendimage{line width=1pt, color=color2}\addlegendentry{Agent 2}
    \addlegendimage{line width=1pt, color=color3}\addlegendentry{Agent 3}
    \addlegendimage{line width=1pt, color=color4}\addlegendentry{Agent 4}
    \addlegendimage{line width=1pt, color=color5}\addlegendentry{Agent 5}
    \addlegendimage{line width=1pt, color=color6}\addlegendentry{Agent 6}
  \end{axis}
\end{tikzpicture}

\maketitle
\begin{abstract}
Optimizing the trajectories of multiple quadrotors in a shared space is a core challenge in various applications. Many existing trajectory optimization methods enforce constraints only at the discretization points, leading to violations between discretization points. They also often lack warm-starting strategies for iterative solution methods such as sequential convex programming, causing slow convergence or sensitivity to the initial guess. We propose a framework for optimizing multiagent quadrotor trajectories that combines a sequential convex programming approach with filtering-based warm-starting. This framework not only ensures constraint satisfaction along the entire continuous-time trajectory but also provides an online warm-starting strategy that accelerates convergence and improves solution quality in numerical experiments. The key idea is to first transform continuous-time constraints into auxiliary nonlinear dynamics and boundary constraints, both of which are compatible with sequential convex programming. Furthermore, we propose a novel warm-starting strategy by approximating the trajectory optimization problem as a Bayesian state estimation problem. This approximation provides an efficient estimate of the optimal trajectories. We demonstrate the proposed framework on multiagent quadrotor trajectory optimization with collision avoidance constraints. Compared with benchmark methods, the proposed framework achieves not only continuous-time constraint satisfaction, but also reduces computation time by up to two orders of magnitude.

\end{abstract}

\section*{Nomenclature}


{\renewcommand\arraystretch{1.0}
\noindent\begin{longtable*}{@{}l @{\quad=\quad} l@{}}
$\bm{c}_{l}$ & the center of no-fly zone for \( l \)-th obstacle ,\,$\mathrm{m}$ \\
$d$   & the minimum allowable distance between any two agents,\,$\mathrm{m}$ \\
$m$   & the total number of quadcopters,\,unitless \\
$m_q$   & quadcopter mass,\,$\mathrm{kg}$ \\
$\bm{n}$ &  a unit vector in the vertical direction:$\begin{bmatrix} 0 & 0 & 1 \end{bmatrix}^\top$,\, unitless \\
$n_o$  & the total number of obstacles,\,unitless \\
$ \bm{r}_{\max}$ & the upper bound vector for the position,\,$\mathrm{m}$ \\
$ \bm{r}_{\min}$ & the lower bound vector for the position,\,$\mathrm{m}$ \\
$T_{\text{max}}$& the maximal thrust magnitude,\,$\mathrm{N}$ \\
$T_{\text{min}}$& the minimal thrust magnitude,\,$\mathrm{N}$ \\
$\bm{u}_{\text{max}}$& the maximal thrust rate vector,\,$\mathrm{N/s}$ \\
$\bm{u}_{\text{min}}$& the minimal thrust rate vector,\,$\mathrm{N/s}$ \\
$v_{\text{max}}$& the maximal velocity magnitude,\,$\mathrm{m/s}$ \\
$\alpha_1$ &the weighting parameter on final time,\, $\mathrm{1/s}$ \\
$\alpha_2$ &the weighting parameter on thrust rate,\,$\mathrm{\frac{s}{N^2}}$ \\

$\alpha_3$ &the weighting parameter on thrust,\, $\mathrm{\frac{1}{s\,N^2}}$ \\


$\theta_{\text{max}}$ & maximum tilting angle,\, $\mathrm{rad}$\\
$\rho_{l}$ & the no-fly zone radius for \( l \)-th obstacle ,\,$\mathrm{m}$ \\

\multicolumn{2}{@{}l}{Matrices}\\
$I_n$ & $n\times n$ identity matrix\\
$0_{m\times n}$ & $m\times n$ zero matrix\\
$\bm{0}_{m}$ & $m\times m$ zero matrix\\
\multicolumn{2}{@{}l}{Vectors}\\
$0_n$ & $n\times 1$ zero vector\\
$\ones_{n}$ & $n\times 1$ one vector\\
\end{longtable*}}

\section{Introduction}
Optimizing the trajectories of multiple quadrotors in cluttered spaces is a core challenge in various applications, such as package delivery and cooperative surveillance~\cite{marcucci2023motion,honig2018trajectory}. These optimal trajectories must minimize cost functions---such as mission time and battery usage---while satisfying various physical and operational constraints, including dynamics and collision avoidance constraints~\cite{huang2012survey}. A principled mathematical optimization framework enables the efficient generation of safe, feasible, and cost-effective trajectories, which is particularly important in multiagent quadrotor systems~\cite{huang2012survey,choi2019multi,wang2024survey,malyuta2022convex}.


Existing approaches to multiagent trajectory optimization vary in their choice of variables for modeling collision avoidance. For example, some approaches formulate collision avoidance with continuous variables and solve the problem via sequential convex programming, which approximates the nonconvex problem as a sequence of convex subproblems~\cite{augugliaro2012generation,XU2022664}. Other approaches use deep reinforcement learning to train a value network with continuous parameters that generates collision‑free trajectories~\cite{7989037}. In contrast, some approaches use mixed-integer programming---a class of optimization problems that includes both continuous and discrete variables---where discrete integer variables encode collision avoidance constraints or agent‑flow assignments~\cite{schouwenaars2001mixed,choi2019multi,nguyen2022collision}. Some variations of mixed-integer programming approaches partition trajectory optimization into two phases: integer programming determines agent‑flow assignment, followed by continuous optimization that smooths the trajectories and enforces dynamic feasibility~\cite{honig2018trajectory}.




One limitation of many existing methods is the risk of constraint violations between discretized trajectory points, also known as \emph{inter-sample constraint violations}. These methods~\cite{augugliaro2012generation,XU2022664,schouwenaars2001mixed,nguyen2022collision,honig2018trajectory,choi2019multi,7989037} discretize a continuous-time trajectory into multiple points and enforce constraints---such as collision avoidance---only at those discretization points, which can lead to violations between consecutive discretization points. These violations are particularly problematic in multiagent trajectory optimization where multiple pairwise agent-agent collision avoidance constraints are present. While increasing the number of discretization points can reduce the risk of such violations, it also increases the problem’s dimensionality by introducing more variables in trajectory optimization, leading to higher computational cost. In mixed-integer programming~\cite{schouwenaars2001mixed,choi2019multi,nguyen2022collision}, this cost increase is further aggravated by the need to introduce additional integer variables at each discretization point.




Another limitation of existing methods in multiagent quadrotor trajectory optimization problems is the lack of efficient warm‐starting strategies for iterative solution methods. These optimization problems are typically nonconvex (due to, \eg, collision avoidance constraints). Most existing iterative solution methods---such as interior‑point methods~\cite{el1996formulation,byrd1999interior,lukvsan2004interior}, sequential quadratic programming~\cite{boggs1995sequential,gill2011sequential}, and sequential convex programming~\cite{morgan2014model,augugliaro2012generation,XU2022664}---compute locally optimal solutions by solving a sequence of subproblems. Each subproblem approximates the original nonconvex problem based on the solution from the previous subproblem. Hence, the quality of the entire sequence of subproblems depends on the first one, which in turn is often determined by an initial guess of the optimal trajectory. Without an initial trajectory sufficiently close to the optimum---also known as a warm-starting trajectory---these methods often suffer from slow convergence or converge to suboptimal trajectories.


Most existing warm-starting strategies for trajectory optimization require costly offline computation. For example, learning-based warm-starting strategies---such as neural-network predictors~\cite{BanerjeeLearning-basedWarm-Starting} or motion memory methods~\cite{lembono2020memory}---provide warm‐starting trajectories by training on large datasets of previously solved problems. However, they often require solving a large number of example problems to construct the training dataset, along with extensive computation to train the model parameters. These demands make them difficult to scale as the number of agents increases in multiagent systems.



We propose a novel framework for multiagent quadrotor systems that achieves two goals: (i) ensuring continuous-time constraint satisfaction, and (ii) including an online, training-free warm-starting strategy that accelerates convergence and improves solution quality in numerical experiments. To achieve the first goal, we propose a sequential convex programming approach that enables continuous-time constraint satisfaction. It transforms continuous-time constraints into auxiliary nonlinear dynamics and final state constraints. Based on this transformation, we propose to compute the optimal trajectories by solving a nonlinear program via the prox-linear method, a first-order method with global convergence guarantees under mild conditions~\cite{drusvyatskiy2018error,drusvyatskiy2019efficiency}.



To achieve the second goal---namely online warm-starting---we propose a filtering‑based technique that extends our prior work~\cite{yuan2024filtering} on discrete-time, fixed-final-time trajectory optimization to the continuous‑time, variable‑final‑time setting. This online warm-starting strategy leverages constraint‑aware particle filtering~\cite{askari2021nonlinear,askari2022sampling,askari2023model}. It first reformulates a trajectory optimization problem as a Bayesian state estimation problem, then estimates the optimal trajectories via efficient nonlinear particle filtering. The estimated trajectories then provide warm-starts for the proposed sequential convex programming approach.

We demonstrate our proposed sequential convex programming approach and filtering-based warm-starting strategy on multiagent quadrotor trajectory optimization subject to collision constraints, including both inter-agent and agent-obstacle collision avoidance. The simulation results show that the trajectories generated by the proposed sequential convex programming approach satisfy all constraints along the entire continuous-time trajectories, not only at the discretization points. Furthermore, compared with benchmark methods, the proposed sequential convex programming approach with filtering-based warm-starting reduces the computation time by up to two orders of magnitude. When compared with random initializations, the filtering-based warm‑starting strategy helps both the proposed sequential convex programming approach and the benchmark methods achieve a lower objective value and a faster reduction in constraint violations.

\section{Continuous-Time Multiagent Quadrotor Trajectory Optimization}

We introduce a continuous-time multiagent quadrotor trajectory optimization problem with variable final time. This optimization problem includes continuous-time three-degree-of-freedom (3-DoF) dynamics models for the \(m \in \mathbb{N}\) quadrotors, path constraints including agent-obstacle and inter-agent collision avoidance constraints, boundary constraints such as initial and final state constraints, and an objective function that optimizes the final time, thrust magnitudes, and thrust rates.


\subsection{Dynamics} \label{subsection: Dynamics}


We consider a 3-DoF continuous-time dynamics model for each quadrotor. In particular, at time \(t \in \mathbb{R}_+\), we let \(\bm{r}^i(t), \bm{v}^i(t) \in \mathbb{R}^3\) denote the position and velocity of the \(i\)-th quadrotor’s center of mass, respectively. We let \(\bm{T}^i(t), \bm{u}^i(t) \in \mathbb{R}^3\) denote the total thrust force and its rate generated by the \(i\)-th quadrotor’s propellers, respectively. Additionally, we let \(m_q\) denote the quadrotor mass. We assume all vehicles are identical and share this value. The vector \( \hat{\bm{T}} = m_q\begin{bmatrix} 0 & 0 & 9.81 \end{bmatrix}^\top \) represents the thrust required to counteract gravity.

We model the continuous-time 3-DoF dynamics of the \(i\)-th quadrotor as follows:
\begin{equation}\label{eqn: quadrotor dynamics}
    \frac{d}{dt} \underbrace{\begin{bmatrix}
        \bm{r}^i(t) \\
        \bm{v}^i(t) \\
        \bm{T}^i(t)
    \end{bmatrix}}_{\bm{x}^i(t)} = \bm{f}(\bm{x}^i(t),\bm{u}^i(t))\coloneqq \begin{bmatrix}
        \bm{0}_{3} & I_3 & \bm{0}_{3} \\
        \bm{0}_{3} &  \bm{0}_{3} & \frac{1}{m_q} I_3 \\
        \bm{0}_{3} &  \bm{0}_{3} & \bm{0}_{3}
    \end{bmatrix}
    \underbrace{\begin{bmatrix}
        \bm{r}^i(t)\\
        \bm{v}^i(t)\\
        \bm{T}^i(t)
    \end{bmatrix}}_{\bm{x}^i(t)} + \begin{bmatrix}
        \bm{0}_{3} \\
        \bm{0}_{3} \\
        I_3
    \end{bmatrix} \bm{u}^i(t) - \begin{bmatrix}
        0_3 \\
        \frac{1}{m_q}\hat{\bm{T}} \\
        0_3
    \end{bmatrix}.
\end{equation}

\subsection{Path Constraints}\label{subsection: Path Constraints}

We consider the following constraints on each quadrotor's position, velocity, thrust and thrust rate:

\begin{enumerate}
\item \textbf{Position.}
The quadrotors must remain within a box
defined by \(\bm{r}_{\min}, \bm{r}_{\max} \in \mathbb{R}^3\), \ie,
\[
  \bm{r}_{\min} \leq \bm{r}^i(t) \leq \bm{r}_{\max},
  \enskip  i = 1,\dots,m.
\]

In addition, we consider the following collision avoidance constraints.

\emph{Obstacle avoidance.} We consider \(n_o\) cylindrical obstacles. We let \(\bm{c}_l \in \mathbb{R}^2\) and \(\rho_l \in \mathbb{R}\) denote the center and radius of the \(l\)-th obstacle, respectively. We introduce the transition matrix \( M_o = \begin{bmatrix}
    I_2 & 0_2
\end{bmatrix}\) to extract the planar position from \(\bm{r}^i(t) \in \mathbb{R}^3\), for all \(i = 1,2,\dots,m\). Based on the above notation, we let
\[
  \rho_l - \norm{M_o\bm{r}^i(t) - \bm{c}_l}_2 \leq 0,
  \enskip  i = 1,\dots,m,
\]
for all \(l=1,2,\dots,n_o\).


\emph{Inter-agent avoidance.}
We let \(d \in \mathbb{R}\) denote the safety distance between any two agents, such that
\[
  d - \norm{\bm{r}^i(t) - \bm{r}^j(t)}_2 \leq 0,
  \enskip  1 \leq i < j \leq m.
\]

\item \textbf{Velocity.}
We let \( v_{\max} \in \mathbb{R} \) denote the quadrotor's maximum velocity magnitude, such that \[ \norm{\bm{v}^i(t)}_2 - v_{\max} \leq 0, \enskip  i = 1,\dots,m. \]

\item \textbf{Thrust.}
The thrust is subject to magnitude and direction constraints. Specifically, the Euclidean norm of the thrust vector is upper bounded by \( T_{\text{max}} \in \mathbb{R} \) and lower bounded by \( T_{\text{min}} \in \mathbb{R} \), such that \[ T_{\text{min}} \leq \|\bm{T}^i(t)\|_2 \leq T_{\text{max}}, \enskip i = 1,\dots,m.  \] 

We consider a tilting angle constraint on the thrust vector (relative to the vertical direction) as an approximation of the quadrotor's attitude limits, without explicitly modeling full attitude dynamics~\cite{yu2023real}. To this end, we define the upper bound of the tilting angle as $\theta_{\max} \in \left[0, \frac{\pi}{2}\right]$, and denote the vertical unit vector by
\[
\bm n = \begin{bmatrix}0 &0 &1\end{bmatrix}^\top.
\]
We formulate the tilting angle constraint as
\[
\cos(\theta_{\max})\cdot \norm{\bm T^i(t)}_2 - \bm n^\top \cdot \bm T^i(t) \leq 0, \enskip i = 1,\dots,m.
\]

\item \textbf{Thrust Rate.}
We bound the thrust rate by \(\bm{u}_{\min},\bm{u}_{\max} \in \mathbb{R}^3\), such that
\[ \bm{u}_{\min} \leq \bm{u}^{i}(t) \leq \bm{u}_{\max}, \enskip i = 1,\dots,m.\]
This ensures that the thrust does not experience extreme changes, thus keeping the trajectory smoother.

\end{enumerate}

\subsection{Boundary Constraints}\label{subsection: Boundary Constraints}

We let \( \hat{\bm x}_0^i,\, \hat{\bm x}_f^i \in \mathbb{R}^9 \) denote the initial and final states for the \(i\)-th quadrotor, respectively, such that
\[
\bm x^i(0) = \hat{\bm x}_0^i,\quad
\bm x^i(t_f) = \hat{\bm x}_f^i,
\quad i = 1,\dots,m.
\]
We also bound the final time by \(t_{\min}, t_{\max}\in \mathbb{R}\), such that
\[
t_{\min}\le t_f\le t_{\max}.
\]

\subsection{Trajectory Optimization}
We introduce the following continuous-time multiagent quadrotor trajectory optimization problem. 
\begin{equation}\label{ocp: min energy and time}
    \begin{array}{ll}
     \underset{\bm{x}^{1:m}, \bm{u}^{1:m},t_f}{\mbox{minimize}} &  \sum_{i=1}^m  \int_0^{t_f} \left(\alpha_1+\alpha_2 \norm{\bm{u}^{i}(t)}_2^2  +  \alpha_3 \norm{\bm{T}^{i}(t)}_2^2\right) \, dt  \\
  \mbox{subject to} & \frac{d}{dt} \bm{x}^{i}(t)=\bm{f}^{i}(\bm{x}^{i}(t), \bm{u}^{i}(t)),\enskip \text{a.e.}~t\in[0, t_f],\enskip i=1, 2, \ldots, m,\\
    &   \bm{r}_{\min} \leq \bm{r}^i(t) \leq \bm{r}_{\max} ,\enskip \text{a.e.}~t\in[0, t_f], \enskip i=1,2, \ldots, m, \\
    &   \norm{\bm{v}^i(t)}_2 \leq v_{\max} ,\enskip \text{a.e.}~t\in[0, t_f], \enskip i=1,2, \ldots, m, \\
    &  T_{\min} \leq \norm{\bm{T}^i(t)}_2 \leq T_{\max} ,\enskip \text{a.e.}~t\in[0, t_f], \enskip i=1,2, \ldots, m, \\
    &  \cos(\theta_{\max})\cdot \norm{\bm{T}^i(t)}_2 - \bm{n}^\top \cdot \bm{T}^i(t) \leq 0 ,\enskip \text{a.e.}~t\in[0, t_f], \enskip i=1,2, \ldots, m, \\
    & \rho_{l} - \norm{M_o \bm{r}^i(t) - \bm{c}_{l}}_2 \leq 0,  \enskip \text{a.e.}~t\in[0, t_f], \enskip l = 1,2, \dots, n_o, \enskip i = 1,2,\dots,m, \\
    & d - \norm{\bm{r}^i(t) - \bm{r}^j(t)}_2 \leq 0, \enskip \text{a.e.}~t\in[0, t_f], \enskip 1 \leq i < j \leq m,\\
    & \bm{u}_{\min}\leq \bm{u}^{i}(t) \leq \bm{u}_{\max}, \enskip \text{a.e.}~t\in[0, t_f], \enskip i=1,2, \ldots, m,\\
   & t_{\min} \leq t_f \leq t_{\max},\\
   & \bm{x}^{i}(0) = \hat{\bm{x}}^{i}_0, \enskip \bm{x}^{i}(t_f)=\hat{\bm{x}}^{i}_f, \enskip i = 1,2, \ldots,m,
    \end{array}
\end{equation}
where \(\bm{x}^{1:m} \coloneqq \{\bm{x}^i\}_{i = 1}^m\), \(\bm{u}^{1:m} \coloneqq \{\bm{u}^i\}_{i = 1}^m\), and \(\bm{x}^i:\mathbb{R}_{\geq 0}\to\mathbb{R}^9\) and  \(\bm{u}^i:\mathbb{R}_{\geq 0}\to\mathbb{R}^3\) are functions of time. We will use the same notation for throughout the paper. The notation a.e.~\(t \in [0, t_f]\) denotes “almost everywhere” in \([0, t_f]\)~\cite{ELANGO2025112464}. The idea of the objective function is to achieve a trade‑off among minimizing the final time, the magnitude of thrust rate over time, and the magnitude of thrust over time. The constraints come from each agent’s dynamics, its path constraints, and its boundary constraints introduced in Section~\ref{subsection: Dynamics}, \ref{subsection: Path Constraints}, and \ref{subsection: Boundary Constraints}. To get weighting parameters, we define:
\[
\alpha_{1}= \tilde{\alpha}_{1}\,\frac{1}{t_{\max}},\quad
\alpha_{2}= \tilde{\alpha}_{2}\,\frac{1}{t_{\max}\norm{\bm{u}_{\max}}_2^2},\quad
\alpha_{3}= \tilde{\alpha}_{3}\,\frac{1}{t_{\max}T_{\max}^{2}},
\]
where each $\tilde{\alpha}_i$ is a normalized weighting parameter and \(\sum_{i=1}^{3}\tilde{\alpha}_i = 1.\)

\section{Sequential Convex Programming with Continuous-Time Constraints}\label{sec: proposed SCP}

We propose a sequential convex programming (SCP) approach to solve the multiagent quadrotor trajectory optimization problem~\eqref{ocp: min energy and time}. This approach ensures constraint satisfaction along the entire continuous-time trajectories. The key idea of this approach is to first transform the continuous-time constraints into auxiliary system dynamics and final state constraints, which lead to a nonlinear program. We then solve this nonlinear program through a sequence of convex optimization problems.


\subsection{Constraint Transformation}
One challenge in solving optimization~\eqref{ocp: min energy and time} is ensuring the inequality constraints on state trajectories. Unlike the input trajectory, the state trajectory depends not only on the system inputs---which correspond to directly controllable actuation---but also on the dynamics equations in~\eqref{eqn: quadrotor dynamics}. To address this challenge, we propose the following transformation that embeds state inequality constraints as auxiliary dynamics and boundary constraints. To this end, at time \(t \in \mathbb{R}_{+}\), we let \(
\bm g^i\left(\bm x^i(t)\right)\colon \mathbb R^9 \to \mathbb R^{10 + n_o}\) for all
\( i = 1,2,\dots,m,
\) denote the function that encodes the inequality constraints for the \(i\)-th agent’s state and its collision avoidance with all obstacles, \ie,
\[
\bm{g}^i\left(\bm{x}^{i}(t)\right) \coloneqq \begin{bmatrix}
    \bm{r}^i(t) - \bm{r}_{\max} \\
    -\bm{r}^i(t) + \bm{r}_{\min} \\
    \norm{\bm{v}^i(t)}_2 - v_{\max} \\
    \norm{\bm{T}^i(t)}_2 - T_{\max} \\
    - \norm{\bm{T}^i(t)}_2 +  T_{\min} \\
    \cos(\theta_{\max})\cdot \norm{\bm{T}^i(t)}_2 - \bm{n}^\top \cdot \bm{T}^i(t) \\
    \rho_{1} - \norm{M_o\bm{r}^i(t) - \bm{c}_{1}}_2 \\
    \vdots \\
    \rho_{n_o} - \norm{M_o\bm{r}^i(t) - \bm{c}_{n_o}}_2 
\end{bmatrix}.
\]
We let \( \bm{h}\left(\bm{r}^{1:m}(t)\right)\colon \mathbb{R}^{3m} \to \mathbb{R}^{\frac{m(m-1)}{2}} \) denote the function that encodes all inter‑agent collision‑avoidance constraints at time \(t \in \mathbb{R}_+\), \ie,
\[
\bm{h}\left(\bm{r}^{1:m}(t)\right) := \begin{bmatrix}
    d - \norm{\bm{r}^1(t) - \bm{r}^2(t)}_2\\
    d - \norm{\bm{r}^1(t) - \bm{r}^3(t)}_2\\
    \vdots\\
    d - \norm{\bm{r}^{m-1}(t) - \bm{r}^m(t)}_2\\
\end{bmatrix},
\]
where \(\bm{r}^{1:m}(t) \coloneqq \{\bm{r}^i(t)\}_{i=1}^m\). Now, we can collect all state‐related inequality constraints in optimization~\eqref{ocp: min energy and time} into \(\bm{g}\left(\bm{x}^{1:m}(t)\right):\mathbb{R}^{9m}\to\mathbb{R}^{n_g}\) at time \(t \in \mathbb{R}_+\), \ie,
\[
 \bm{g}\left(\bm{x}^{1:m}(t)\right) := \begin{bmatrix}
     \bm{g}^1\left(\bm{x}^1(t)\right)\\
     \vdots\\
     \bm{g}^m\left(\bm{x}^m(t)\right)\\
     \bm{h}\left(\bm{r}^{1:m}(t)\right)
 \end{bmatrix}.
\]

To transform the continuous-time inequality constraints \(\bm{g}\left(\bm{x}^{1:m}(t)\right)\) into dynamics constraints, we introduce the following auxiliary states, such that 
\begin{subequations}
    \begin{align}
     &y(t)= \int_0^t \norm{\max\{\bm{g}(\bm{x}^{1:m}(\sigma)), 0_{n_g}\}}_2^2 \, d\sigma,\\
       & w(t)=\sum_{i=1}^m \int_0^t \left( \alpha_1 + \alpha_2 \norm{\bm{u}^{i}(\sigma)}_2^2 + \alpha_3 \norm{\bm{T}^{i}(\sigma)}_2^2 \right)\,d\sigma,
    \end{align}
\end{subequations}
where \(y(t)\) represents the integrated constraint violation, and \(w(t)\) corresponds to the objective function when \(t = t_f\). By enforcing that the integrated constraint violation satisfies \(y(t_f) = 0\), we guarantee that constraint satisfaction holds almost everywhere over time. Thus, we reformulate the original optimization problem~\eqref{ocp: min energy and time} as
\begin{equation}\label{ocp: min_energy_Mayer}
    \begin{array}{ll}
     \underset{\bm{x}^{1:m}, \bm{u}^{1:m},w,y,t_f}{\mbox{minimize}} &  w(t_f) \\
  \mbox{subject to} & \frac{d}{dt} \begin{bmatrix}
  \bm{x}^{1}(t)\\
  \vdots\\
  \bm{x}^{m}(t)\\
  y(t)\\
  w(t)
  \end{bmatrix}=\begin{bmatrix}
      \bm{f}^{1}(\bm{x}^{1}(t), \bm{u}^{1}(t))\\
      \vdots\\
      \bm{f}^{m}(\bm{x}^{m}(t), \bm{u}^{m}(t))\\
       \norm{\max\{\bm{g}(\bm{x}^{1:m}(t)), 0_{n_g}\}}_2^2\\
       m\alpha_1 + \alpha_2 \sum_{i=1}^m \norm{\bm{u}^{i}(t)}_2^2+\alpha_3 \sum_{i=1}^m \norm{\bm{T}^i(t)}_2^2
  \end{bmatrix},\enskip \text{a.e.}~t\in[0, t_f],\\
   & \bm{u}_{\min} \leq \bm{u}^{i}(t) \leq \bm{u}_{\max}, \enskip \text{a.e.}~t\in[0, t_f],\enskip i=1,2, \ldots, m,\\
   & t_{\min} \leq t_f \leq t_{\max},\\
    & \bm{x}^{i}(0) = \hat{\bm{x}}^{i}_0,\enskip
      \bm{x}^{i}(t_f)=\hat{\bm{x}}^{i}_f,\enskip i=1,2, \ldots, m,\\
      & y(0)=0,\enskip w(0)=0,\enskip y(t_f)=0. 
    \end{array}
\end{equation}


\subsection{From Variable-Final-Time to Fixed-Final-Time}
To facilitate numerical optimization, we need to parameterize the trajectory in optimization~\eqref{ocp: min_energy_Mayer}---which has a variable final time---using a fixed number of parameters. To this end, we first transform it into an equivalent fixed-final-time problem, then discretize it into a finite number of steps. In particular, we apply the \emph{time-scaling} transformation~\cite{Szmuk,ELANGO2025112464,KAMATH20233118} by introducing a new variable \(\tau\in[0, 1]\), named \emph{time-scaling control}. We relate this new variable to the original time variable \(t\in[0, t_f]\) through the following boundary value problem: 
\begin{equation}
\begin{aligned}
    \frac{d t(\tau)}{d\tau} & = s(\tau), \enskip \tau\in[0, 1],\\
    t(0)& =0, \enskip t(1)=t_f.
\end{aligned}
\end{equation}
In addition, we treat \(s(\tau)\) as a continuous-time control input and introduce the following notion of \emph{augmented state and input variables} 
\begin{equation}
\begin{aligned}
    & \overline{\bm{x}}(\tau)\coloneqq \begin{bmatrix}
         \bm{x}^1(t(\tau)) ^\top & \dots &  \bm{x}^{m}(t(\tau)) ^\top & y(t(\tau)) & w(t(\tau))
    \end{bmatrix}^\top \in \mathbb{R}^{n_x+2} , \\
    & \overline{\bm{u}}(\tau)\coloneqq \begin{bmatrix}
        \bm{u}^1(t(\tau))^\top & \dots & \bm{u}^{m}(t(\tau))^\top &s(\tau)
    \end{bmatrix}^\top \in \mathbb{R}^{n_u+1},
\end{aligned}
\end{equation}
where \( n_x = 9m \) and \( n_u = 3m \) represent the total number of states and inputs of all \(m\) agents, respectively. By applying the chain rule, we can show the following
\begin{equation}
\begin{aligned}
    \frac{d}{d\tau} \overline{\bm{x}}(\tau) & = \begin{bmatrix}
         \frac{d}{d\tau}\bm{x}^1(t(\tau))\\
         \vdots\\
        \frac{d}{d\tau}\bm{x}^{m}(t(\tau))\\
        \frac{d}{d\tau} y(t(\tau))\\
        \frac{d}{d\tau} w(t(\tau)) 
    \end{bmatrix}= \begin{bmatrix}
         \frac{d}{dt(\tau)}\bm{x}^1(t(\tau))\frac{d t(\tau)}{d\tau}\\
         \vdots\\
         \frac{d}{dt(\tau)}\bm{x}^m(t(\tau))\frac{d t(\tau)}{d\tau}\\
        \frac{d}{dt(\tau)} y(t(\tau))\frac{d t(\tau)}{d\tau}\\
        \frac{d}{dt(\tau)} w(t(\tau))\frac{d t(\tau)}{d\tau} \end{bmatrix}
        =s(\tau)\begin{bmatrix}   
          \bm{f}^1(\bm{x}^1(t(\tau)), \bm{u}^1(t(\tau)))\\
          \vdots\\
          \bm{f}^m(\bm{x}^m(t(\tau)),\bm{u}^m(t(\tau)))\\
       \norm{\max\{\bm{g}(\bm{x}^{1:m}(t(\tau))), 0_{n_g}\}}_2^2\\
       m\alpha_1 + \alpha_2 \sum_{i=1}^m \norm{\bm{u}^i(t(\tau))}_2^2+ \alpha_3 \sum_{i=1}^m \norm{\bm{T}^{i}(t(\tau))}_2^2\end{bmatrix}\\
       &\coloneqq \overline{\bm{f}}(\overline{\bm{x}}(\tau), \overline{\bm{u}}(\tau)).           
\end{aligned}
\label{eq: augmented dynamics}
\end{equation}
Furthermore, one can verify that if
\begin{equation}\label{eq: timing-factor constraint}
    t_{\min}\leq s(\tau)\leq t_{\max}, \enskip \tau\in[0, 1],
\end{equation}
then 
\begin{equation}
    t_{\min}\leq t_f=t(1)=t(0)+\int_0^1 s(\tau) \, d\tau\leq t_{\max},
\end{equation}
where the two inequalities become equalities if \(s(\tau)=t_{\min}\) and \(s(\tau)=t_{\max}\) for all \(\tau\in[0, 1]\), respectively. This continuous-time inequality~\eqref{eq: timing-factor constraint} enables the adaptive time discretization scheme in Section~\ref{subsec: Discretization}.


Based on the augmented dynamics in~\eqref{eq: augmented dynamics}, we can rewrite optimization~\eqref{ocp: min_energy_Mayer} equivalently as follows:
\begin{equation}\label{opt: time-scaled}
    \begin{array}{ll}
     \underset{\overline{\bm{x}}, \overline{\bm{u}}}{\mbox{minimize}} & \begin{bmatrix}
        0_{n_x +1}^\top & 1
    \end{bmatrix} \overline{\bm{x}}(1)\\
  \mbox{subject to} & \frac{d}{d\tau} \overline{\bm{x}}(\tau)=\overline{\bm{f}}(\overline{\bm{x}}(\tau), \overline{\bm{u}}(\tau)),\enskip \text{a.e.}~\tau\in[0, 1],\\
  & \overline{\bm{x}}(0)= \overline{\bm{x}}_0, \enskip \begin{bmatrix}
        I_{n_x+1} & 0_{n_x+1}
    \end{bmatrix} \left(\overline{\bm{x}}(1) - \overline{\bm{x}}_f\right)= 0_{n_x+1},\\
  & \overline{\bm{u}}_{\min}\leq \overline{\bm{u}}(\tau) \leq \overline{\bm{u}}_{\max}, \enskip \text{a.e.}~\tau\in[0, 1],
    \end{array}
\end{equation}
where
\begin{align*}
&\overline{\bm{x}}_0 = \begin{bmatrix}
    \left(\hat{\bm{x}}^{1}_0\right)^\top & \dots &  \left(\hat{\bm{x}}^{m}_0\right)^\top & 0_2^\top
\end{bmatrix}^\top, \enskip  \overline{\bm{x}}_f = \begin{bmatrix}
        \left(\hat{\bm{x}}^{1}_f\right)^\top & \dots & \left(\hat{\bm{x}}^{m}_f\right)^\top & 0_2^\top \end{bmatrix}^\top,\\
&\overline{\bm{u}}_{\min} = \begin{bmatrix}
        \left( \ones_m \otimes \bm{u}_{\min}\right)^\top &  t_{\min}
    \end{bmatrix}^\top, \enskip \overline{\bm{u}}_{\max} = \begin{bmatrix}
    \left( \ones_m \otimes \bm{u}_{\max} \right)^\top & t_{\max}
    \end{bmatrix}^\top,
\end{align*}
and \(\otimes\) is the Kronecker product. Notice that changing the value of the last coordinate in \(\overline{\bm{x}}_f\) does not affect the final state constraints in optimization~\eqref{opt: time-scaled}. Here we choose the value zero for simplicity.

\subsection{Parameterization and Discretization}\label{subsec: Discretization}

We now parameterize optimization~\eqref{opt: time-scaled} as a finitely-dimensional nonlinear program. To this end, we parameterize the input as a piecewise constant function of \(\tau\). Let \(\tau_1, \tau_2,\ldots, \tau_N\in\mathbb{R}\) be such that
\begin{equation}
    \tau_k = \frac{k-1}{N-1},\enskip k = 1,2,\ldots,N.
\end{equation}
In addition, let \(\bm{\eta}_1, \bm{\eta}_2, \ldots, \bm{\eta}_{N-1}\in\mathbb{R}^{n_u+1}\) denote the \emph{input parameters}, such that
\begin{equation}
    \overline{\bm{u}}(\tau)=\bm{\eta}_k, \enskip \tau\in[\tau_k, \tau_{k+1}),
\end{equation}
for all \(k=1, 2, \ldots, N-1\). Next, let \(\bm{\xi}_1, \bm{\xi}_2, \ldots, \bm{\xi}_N\in\mathbb{R}^{n_x+2}\) denote the \emph{state parameters}. For any \(k=1, 2, \ldots, N-1\), we define function \(\bm{F}_k:\mathbb{R}^{n_x+2}\times \mathbb{R}^{n_u+1}\to\mathbb{R}^{n_x+2}\) such that
\begin{equation}
\begin{aligned}
    \bm{\xi_{k+1}} = \overline{\bm{x}}(\tau_{k+1}) = \bm{F}_k(\bm{\xi}_k, \bm{\eta}_k)  
    &= \overline{\bm{x}}(\tau_{k}) + \int_{\tau_k}^{\tau_{k+1}} \overline{\bm{f}}(\overline{\bm{x}}(\tau), \bm{\eta}_k) \, d\tau = \bm{\xi_k} + \int_{\tau_k}^{\tau_{k+1}} \overline{\bm{f}}(\overline{\bm{x}}(\tau), \bm{\eta}_k) \, d\tau.
\end{aligned}
\end{equation}
With these parameters and definitions, we can rewrite optimization~\eqref{opt: time-scaled} equivalently as follows
\begin{equation}\label{opt: min time}
    \begin{array}{ll}
    \underset{\bm{\xi}_{1:N}, \bm{\eta}_{1:N-1}}{\mbox{minimize}} & \begin{bmatrix}
        0_{n_x+1}^\top & 1 
    \end{bmatrix} \bm{\xi}_N \\
    \mbox{subject to}  & \bm{\xi}_{k+1}=\bm{F}_k(\bm{\xi}_k, \bm{\eta}_k), \enskip k=1, 2, \ldots, N-1,\\
    &\bm{\xi}_1=\overline{\bm{x}}_0, \enskip \begin{bmatrix}
        I_{n_x+1} & 0_{n_x+1}
    \end{bmatrix} \left(\bm{\xi}_N - \overline{\bm{x}}_f \right)= 0_{n_x+1},\\
    & \overline{\bm{u}}_{\min} \leq \bm{\bm{\eta}}_k\leq \overline{\bm{u}}_{\max}, \enskip k=1, 2, \ldots, N-1,
    \end{array}
\end{equation}
where \(\bm{\xi}_{1:N} \coloneqq \{\bm{\xi}_k\}_{k=1}^N\) and \(\bm{\eta}_{1:N-1} \coloneqq \{\bm{\eta}_k\}_{k=1}^{N-1}\). We will use the same notation throughout the paper.

The final state constraint in optimization~\eqref{opt: min time} can lead to a violation of the linear independence constraint qualification~\cite[Sec.~3.1]{ELANGO2025112464}. In practice, we often relax the final state constraint in \eqref{opt: min time} and solve the following optimization problem instead:
\begin{equation}\label{opt: min_time_relaxed_v1}
    \begin{array}{ll}
    \underset{\bm{\xi}_{1:N}, \bm{\eta}_{1:N-1}}{\mbox{minimize}} & \begin{bmatrix}
        0_{n_x+1}^\top & 1 
    \end{bmatrix} \bm{\xi}_N \\
    \mbox{subject to}  & \bm{\xi}_{k+1}=\bm{F}_k(\bm{\xi}_k, \bm{\eta}_k), \enskip k=1, 2, \ldots, N-1,\\
    &\bm{\xi}_1=\overline{\bm{x}}_0, \enskip \begin{bmatrix}
        I_{n_x} & 0_{n_x\times 2}
    \end{bmatrix} \left(\bm{\xi}_N - \overline{\bm{x}}_f \right) = 0_{n_x},\\
    & \begin{bmatrix}
        0_{n_x}^\top & 1 & 0
    \end{bmatrix}\bm{\xi}_k \leq \gamma, \enskip k = 2,\ldots,N,\\
    & \overline{\bm{u}}_{\min} \leq \bm{\eta}_k\leq \overline{\bm{u}}_{\max}, \enskip k=1, 2, \ldots, N-1,
    \end{array}
\end{equation}
where \(\gamma\in\mathbb{R}_{>0}\) is a relaxation parameter (\eg, \(\gamma=10^{-6}\)). 

Due to the exact penalty theorem~\cite[Theorem.~17.4]{nocedal1999numerical}, the optimization in~\eqref{opt: min_time_relaxed_v1} is equivalent to the following one for some sufficiently large \(\beta\in\mathbb{R}_{>0}\): 
\begin{equation}\label{opt: min_time_relaxed_prox}
    \begin{array}{ll}
    \underset{\substack{\bm{\xi}_{1:N}, \bm{\eta}_{1:N-1},\\
    \bm{q}_{1:N-1},\bm{z}_{1:N-1}}}{\mbox{minimize}} & \begin{bmatrix}
        0_{n_x+1}^\top & 1 
    \end{bmatrix} \bm{\xi}_N + \beta \sum_{k=1}^{N-1} \left( \ones_{n_x+2}^\top (\bm{q}_k+\bm{z}_k)\right)  \\
    \mbox{subject to}  & \bm{\xi}_{k+1}= \bm{F}_k(\bm{\xi}_k, \bm{\eta}_k) + \bm{q}_k - \bm{z}_k, \enskip k=1, 2, \ldots, N-1,\\
    &\bm{\xi}_1 = \overline{\bm{x}}_0, \enskip \begin{bmatrix}
        I_{n_x} & 0_{n_x\times 2}
    \end{bmatrix} \left(\bm{\xi}_N - \overline{\bm{x}}_f \right)= 0_{n_x},\\
    & \begin{bmatrix}
        0_{n_x}^\top & 1 & 0
    \end{bmatrix}\bm{\xi}_k \leq \gamma,\enskip k = 2,\ldots,N,\\
    & \overline{\bm{u}}_{\min} \leq \bm{\eta}_k\leq \overline{\bm{u}}_{\max}, \enskip k=1, 2, \ldots, N-1,\\
    &\bm{q}_k \geq 0_{n_x+2},\enskip \bm{z}_k \geq 0_{n_x+2},\enskip k=1, 2, \ldots, N-1.
    \end{array}
\end{equation}
Compared with optimization~\eqref{opt: min_time_relaxed_v1}, optimization~\eqref{opt: min_time_relaxed_prox} allows violations of the nonlinear dynamics constraints by introducing two slack variables, \(\bm{q}_k \in \mathbb{R}^{n_x+2}_{\geq 0}\) and \(\bm{z}_k \in \mathbb{R}^{n_x+2}_{\geq 0}\), but penalizes these violations through an additional penalty term in the objective function. The benefit of introducing these two slack variables is to prevent \emph{artificial infeasibility} arising from SCP~\cite{malyuta2022convex}.



\subsection{Linearization and Prox-Linear Method}\label{sec: linearization}

The difficulty in solving optimization~\eqref{opt: min_time_relaxed_v1} comes from the nonlinear dynamics constraints, particularly due to the lack of a closed-form formula for the function \(\bm{F}_k\). To address this challenge, we approximate the optimization~\eqref{opt: min_time_relaxed_v1} by constructing a sequence of convex subproblems---a process known as SCP---where each subproblem computes a linear approximation of the nonlinear dynamics based on the solution from the previous subproblem. In the following, we first show how to compute such linear approximations by solving ordinary differential equations, and then discuss how to construct a sequence of subproblems via the prox-linear method, which enjoys global convergence guarantees under mild regularity conditions~\cite{drusvyatskiy2018error,drusvyatskiy2019efficiency}. 


Given linearization points \(\hat{\bm{\xi}}_k \in \mathbb{R}^{n_x+2}\) and \( \hat{\bm{\eta}}_k \in \mathbb{R}^{n_u+1}\), our goal is to compute the following first-order approximation of \(\bm{F}_k(\bm{\xi}_k, \bm{\eta}_k)\): 
\begin{equation}
    \bm{F}_k(\bm{\xi}_k, \bm{\eta}_k)\approx A_k\bm{\xi}_k+B_k\bm{\eta}_k+\bm{c}_k,
\end{equation}
where
\begin{subequations}\label{eqn: AkBkck}
    \begin{align}
    A_k &\coloneqq \frac{\partial }{\partial \bm{\xi}_k}\bm{F}_k(\bm{\xi}_k, \bm{\eta}_k)\left|\begin{smallmatrix}\bm{\xi}_k=\hat{\bm{\xi}}_k\\
    \bm{\eta}_k=\hat{\bm{\eta}}_k
    \end{smallmatrix}\right. =\frac{\partial}{\partial \bm{\xi}_k} \overline{\bm{x}}(\tau_{k+1})\left|\begin{smallmatrix}
    \bm{\xi}_k=\hat{\bm{\xi}}_k\\
    \bm{\eta}_k=\hat{\bm{\eta}}_k
    \end{smallmatrix}\right.,\\
    B_k &\coloneqq \frac{\partial }{\partial \bm{\eta}_k}\bm{F}_k(\bm{\xi}_k, \bm{\eta}_k)\left|\begin{smallmatrix}\bm{\xi}_k=\hat{\bm{\xi}}_k\\
    \bm{\eta}_k=\hat{\bm{\eta}}_k
    \end{smallmatrix}\right.=\frac{\partial}{\partial \bm{\eta}_k} \overline{\bm{x}}(\tau_{k+1})\left|\begin{smallmatrix}
    \bm{\xi}_k=\hat{\bm{\xi}}_k\\
    \bm{\eta}_k=\hat{\bm{\eta}}_k
    \end{smallmatrix}\right.,\\
    \bm{c}_k&\coloneqq - A_k\hat{\bm{\xi}}_k-B_k\hat{\bm{\eta}}_k + \bm{F}_k(\bm{\xi}_k, \bm{\eta}_k)\left|\begin{smallmatrix}\bm{\xi}_k=\hat{\bm{\xi}}_k\\
    \bm{\eta}_k=\hat{\bm{\eta}}_k
    \end{smallmatrix}\right.=- A_k\hat{\bm{\xi}}_k-B_k\hat{\bm{\eta}}_k+\overline{\bm{x}}(\tau_{k+1})\left|\begin{smallmatrix}
    \bm{\xi}_k=\hat{\bm{\xi}}_k\\
    \bm{\eta}_k=\hat{\bm{\eta}}_k
    \end{smallmatrix}\right..
    \end{align}
\end{subequations}
To this end, for any \(\tau\in[\tau_k,\tau_{k+1})\), let \(\overline{\bm{x}}(\tau)\) be the solution of the following ordinary differential equation:
\begin{subequations}\label{eqn: interval ode}
\begin{align}
    \overline{\bm{x}}(\tau) &= \overline{\bm{x}}(\tau_k)+\int_{\tau_k}^{\tau}\overline{\bm{f}}(\overline{\bm{x}}(\sigma), \overline{\bm{u}}(\sigma)) \, d\sigma,\label{eqn: interval ode1} \\
    \overline{\bm{x}}(\tau_k) & =\bm{\xi}_k,\enskip \overline{\bm{u}}(\tau)=\bm{\eta}_k,\enskip \tau\in[\tau_k, \tau_{k+1}).
\end{align}
\end{subequations}
Next, let \(\Phi_x(\tau)\) and \(\Phi_u(\tau)\) be the Jacobians of \(\overline{\bm{x}}(\tau)\) with respect to \(\bm{\xi}_k\) and \(\bm{\eta}_k\), respectively. 
\begin{equation}
        \Phi_x(\tau) \coloneqq \frac{\partial}{\partial \bm{\xi}_k} \overline{\bm{x}}(\tau), \enskip \Phi_u(\tau) \coloneqq \frac{\partial}{\partial \bm{\eta}_k} \overline{\bm{x}}(\tau).
\end{equation}
The equation in~\eqref{eqn: interval ode} implies that 
\begin{subequations}\label{eqn: Jacobian ode}
    \begin{align}
           \Phi_x(\tau) & =I_{n_x+2}+\int_{\tau_k}^\tau\frac{\partial }{\partial \overline{\bm{x}}(\sigma)}\overline{\bm{f}}(\overline{\bm{x}}(\sigma), \overline{\bm{u}}(\sigma))\frac{\partial}{\partial \bm{\xi}_k} \overline{\bm{x}}(\sigma) \, d\sigma,\\
           \Phi_u(\tau) & =\int_{\tau_k}^\tau\left(\frac{\partial }{\partial \overline{\bm{x}}(\sigma)}\overline{\bm{f}}(\overline{\bm{x}}(\sigma), \overline{\bm{u}}(\sigma))\frac{\partial}{\partial \bm{\eta}_k} \overline{\bm{x}}(\sigma)+\frac{\partial}{\partial \overline{\bm{u}}(\sigma)}\overline{\bm{f}}(\overline{\bm{x}}(\sigma), \overline{\bm{u}}(\sigma))\right) \, d\sigma.
    \end{align}
\end{subequations}
By combining \eqref{eqn: AkBkck}, \eqref{eqn: interval ode}, and \eqref{eqn: Jacobian ode}, we conclude that 
\begin{equation}
    \label{eq: A_kB_kc_k}
    \bm{c}_k=- A_k\hat{\bm{\xi}}_k-B_k\hat{\bm{\eta}}_k+\overline{\bm{x}}(\tau_{k+1})\left|\begin{smallmatrix}
    \bm{\xi}_k=\hat{\bm{\xi}}_k\\
    \bm{\eta}_k=\hat{\bm{\eta}}_k
    \end{smallmatrix}\right., \enskip A_k=\Phi_x(\tau_{k+1})\left|\begin{smallmatrix}
    \bm{\xi}_k=\hat{\bm{\xi}}_k\\
    \bm{\eta}_k=\hat{\bm{\eta}}_k
    \end{smallmatrix}\right., \enskip B_k=\Phi_u(\tau_{k+1})\left|\begin{smallmatrix}
    \bm{\xi}_k=\hat{\bm{\xi}}_k\\
    \bm{\eta}_k=\hat{\bm{\eta}}_k
    \end{smallmatrix}\right.,
\end{equation}
where \(\overline{\bm{x}}(\tau), \Phi_x(\tau), \Phi_u(\tau)\) are the solution of the following systems of ordinary differential equations
\begin{subequations}
    \begin{align}
           \frac{d}{d\tau}\overline{\bm{u}}(\tau)&=0_{n_u+1},\\
           \frac{d}{d\tau}\overline{\bm{x}}(\tau)&=\overline{\bm{f}}(\overline{\bm{\bm{x}}}(\tau), \overline{\bm{u}}(\tau)),\\
           \frac{d}{d\tau} \Phi_x(\tau) & =\frac{\partial}{\partial \overline{\bm{x}}(\tau)}\overline{\bm{f}}(\overline{\bm{x}}(\tau), \overline{\bm{u}}(\tau))\Phi_x(\tau), \\
            \frac{d}{d\tau} \Phi_u(\tau) &=\frac{\partial}{\partial \overline{\bm{x}}(\tau)}\overline{\bm{f}}(\overline{\bm{x}}(\tau), \overline{\bm{u}}(\tau))\Phi_u(\tau)+\frac{\partial }{\partial \overline{\bm{u}}(\tau)}\overline{\bm{f}}(\overline{\bm{x}}(\tau), \overline{\bm{u}}(\tau)), \\
            \overline{\bm{x}}(\tau_k)&=\hat{\bm{\xi}}_k,\enskip \overline{\bm{u}}(\tau_k)=\hat{\bm{\eta}}_k,\enskip \Phi_x(\tau_k)=I_{n_x+2},\enskip \Phi_u(\tau_k)=0_{(n_x+2)\times (n_u+1)}.
    \end{align}
\end{subequations}

The prox‑linear method solves nonconvex problems by iteratively linearizing the nonlinear constraints and enjoys a global convergence guarantee under mild regularity conditions~\cite{drusvyatskiy2018error,drusvyatskiy2019efficiency}. We apply the prox‑linear method to optimization~\eqref{opt: min_time_relaxed_prox} via iterative linearization, where the introduced slack variables help avoid \emph{artificial infeasibility}. Let \( (\{\hat{\bm{\xi}}_k^{(j)}\}_{k=1}^N,\{\hat{\bm{\eta}}_k^{(j)}\}_{k=1}^{N-1} )\) denotes the solution at the \(j\)-th iteration (\(j \in \mathbb{N}\)) of prox-linear method, and \(A_k^{(j)}\), \(B_k^{(j)}\), and \(\bm{c}_k^{(j)}\) is obtained through~\eqref{eq: A_kB_kc_k} by letting \( \hat{\bm{\xi}}_k = \hat{\bm{\xi}}_k^{(j)}\) and \( \hat{\bm{\eta}}_k = \hat{\bm{\eta}}_k^{(j)}\) for all \(k = 1,2,\ldots,N-1\) . At (\(j+1\))-th iteration, the prox-linear method solves the following approximation of optimization~\eqref{opt: min_time_relaxed_prox}:
\begin{equation}\label{opt: prox-linear iter}
    \begin{array}{ll}
    \underset{\substack{\bm{\xi}_{1:N}, \bm{\eta}_{1:N-1},\\
    \bm{q}_{1:N-1},\bm{z}_{1:N-1}}}{\mbox{minimize}} & \begin{bmatrix}
        0_{n_x+1}^\top & 1 
    \end{bmatrix} \bm{\xi}_N + \beta \sum_{k=1}^{N-1} \left( \ones_{n_x+2}^\top (\bm{q}_k+\bm{z}_k)\right) +\frac{1}{2\rho}\sum_{k=1}^N\norm{\bm{\xi}_k-\hat{\bm{\xi}}_k^{(j)}}_2^2 +\frac{1}{2\rho}\sum_{k=1}^{N-1}\norm{\bm{\eta}_k-\hat{\bm{\eta}}_k^{(j)}}_2^2\\
    \quad\mbox{subject to}  & \bm{\xi}_{k+1}=A_k^{(j)}\bm{\xi}_k + B_k^{(j)}\bm{\eta}_k+\bm{c}_k^{(j)}+\bm{q}_k-\bm{z}_k, \enskip k= 1, 2,\ldots, N-1,\\
     & \bm{\xi}_1=\overline{\bm{x}}_0, \enskip \begin{bmatrix}
        I_{n_x} & 0_{n_x\times 2}
    \end{bmatrix}\left(\bm{\xi}_N - \overline{\bm{x}}_f \right)= 0_{n_x},\\
    & \begin{bmatrix}
        0_{n_x}^\top & 1 & 0
    \end{bmatrix}\bm{\xi}_k \leq \gamma, \enskip k = 2,\ldots,N,\\
    & \overline{\bm{u}}_{\min} \leq \bm{\eta}_k\leq \overline{\bm{u}}_{\max}, \enskip k=1, 2, \ldots, N-1,\\
    & \bm{q}_k\geq 0_{n_x+2},\enskip \bm{z}_k\geq 0_{n_x+2},\enskip  k=1, 2, \ldots, N-1,
    \end{array}
\end{equation}
where \(\rho\) is a weighting parameter. We summarize the prox-linear method in Algorithm~\ref{alg: prox-linear}. The idea is to linearize the nonlinear constraints around the nominal trajectory obtained in the previous iteration and to penalize deviations from this nominal trajectory in the objective function.

\begin{algorithm}
    \caption{Prox-Linear Method}
    \label{alg: prox-linear}
    \begin{algorithmic}[1]
        \Require initial trajectory \( (\{\hat{\bm{\xi}}_k^{(0)}\}_{k=1}^N,\{\hat{\bm{\eta}}_k^{(0)}\}_{k=1}^{N-1} )\), \(j=0\) , penalty weight \(\beta \in \mathbb{R}_{>0}\), maximum number of iterations \(j_{\max} \in \mathbb{N}\), tolerance \(\epsilon_\text{tol} \in \mathbb{R}_{>0}\)
        \While{\( j \leq j_{\max}\)}
        \State Let \(\rho \geq 1/\beta\) in optimization~\eqref{opt: prox-linear iter}, solve it and obtain an optimal trajectory \( (\{\bm{\xi}_k^\star\}_{k=1}^N,\{\bm{\eta}_k^\star\}_{k=1}^{N-1}) \)
        \If{\( \sum_{k=1}^{N} \norm{\hat{\bm{\xi}}_k^{(j)} - \bm{\xi}_k^\star}_2^2 + \sum_{k=1}^{N-1} \norm{\hat{\bm{\eta}}_k^{(j)} - \bm{\eta}_k^\star}_2^2 \leq \epsilon_\text{tol} \)}
        \State Break
        \Else
        \State Let \( \{\hat{\bm{\xi}}_k^{(j+1)}\}_{k=1}^N = \{\bm{\xi}_k^\star\}_{k=1}^N\) and \( \{\hat{\bm{\eta}}_k^{(j+1)}\}_{k=1}^{N-1} = \{\bm{\eta}_k^\star\}_{k=1}^{N-1}\)
        \State \(j \gets j+1\)
        \EndIf
        \EndWhile
        \Ensure \( (\{\bm{\xi}_k^\star\}_{k=1}^N,\{\bm{\eta}^\star_k\}_{k=1}^{N-1}) \)
    \end{algorithmic}
\end{algorithm}


\section{Warm-Starting via Constraint-Aware Particle Filtering}\label{sec: proposed warm-starting}

A key challenge in multi‐quadrotor trajectory optimization is the lack of warm-starting trajectories. Iterative solution methods---which solve the nonconvex problem via a sequence of subproblems---often converge slowly or to suboptimal solutions without an initial guess near the optimum. We propose a warm-starting strategy by using constraint-aware particle filtering~\cite{askari2021nonlinear,askari2022sampling,askari2023model} to generate the sampled trajectories that approximate the optimal ones.


\subsection{Control-Estimation Duality}

The key idea of constraint‑aware particle filtering is to formulate an optimization problem as a Bayesian state estimation problem and to estimate the optimal trajectories via efficient particle filtering. However, constraint-aware particle filtering interprets each step’s quadratic tracking error as a Bayesian observation and therefore requires per‑step reference trajectories. To apply constraint‑aware particle filtering to the optimization~\eqref{opt: min_time_relaxed_v1}, we approximate the objective function in optimization~\eqref{opt: min_time_relaxed_v1}---which depends only on the final state---with a quadratic objective function that penalizes deviation from a reference trajectory. This quadratic objective function jointly minimizes thrust rate magnitude, thrust magnitude, final time, deviation from the reference position, and inequality constraint violations. To this end, we introduce the transition matrices \( M_x \in \mathbb{R}^{3m \times (n_x+2)}\) and \( M_T \in \mathbb{R}^{3m \times (n_x+2)} \) map the state \(\bm{\xi}_k\) to column vectors that extract the position and thrust of all agents, respectively. Specifically,
\[
M_x = \begin{bmatrix}
I_m \otimes \begin{bmatrix}
    I_3 & 0_{3\times6}
\end{bmatrix} & 0_{3m \times 2}
\end{bmatrix}, \enskip M_T = \begin{bmatrix}
I_m \otimes \begin{bmatrix}
    0_{3\times6} & I_3
\end{bmatrix} & 0_{3m \times 2}
\end{bmatrix},
\]
where \(\otimes\) is the Kronecker product. Now we define the position reference \(\hat{\bm{x}}_k \in \mathbb{R}^{3m}\) as:
\begin{equation}
\hat{\bm{x}}_k = \frac{N-k}{N-1}  M_x \overline{\bm{x}}_0 + \frac{k-1}{N-1} M_x \overline{\bm{x}}_f,
\end{equation}
for all \(k = 1,2,\dots,N\), to serve as a reference position for all agents at discretization step \(k\). The scalar \( \epsilon \in (0,1) \) is chosen such that the weighting factor on position tracking increases as the position state approaches its desired final state, penalizing errors more heavily near the end of the trajectory. We approximate the optimization~\eqref{opt: min_time_relaxed_v1} as
\begin{equation}\label{pf: MHE_before}
    \begin{array}{ll}
    \underset{\bm{\xi}_{1:N}, \bm{\eta}_{1:N}}{\mbox{minimize}} & \sum_{k=1}^{N} \left(\norm{ \bm{e}_k }_Q^2 + \norm{ [\bm{G}(\bm{\xi}_k,\bm{\eta}_k)]_+ +\nu \bm{1}_{n_{G}} }_2^2 + \norm{\bm{\eta}_k - \begin{bmatrix}
        0_{n_u}^\top & t_{\min}
    \end{bmatrix}^\top}_R^2 \right) \\
    \mbox{subject to} & \bm{\xi}_{k+1} = \bm{F}_k(\bm{\xi}_k, \bm{\eta}_k), \enskip k = 1, 2, \ldots, N-1, \\
    & \bm{\xi}_1 = \overline{\bm{x}}_0, \\
    & \bm{e}_k = \underbrace{\begin{bmatrix} 
       \epsilon^{\frac{N-k}{2}}\hat{\bm{x}}_k \\
        0_{3m}
    \end{bmatrix}}_{\hat{\bm{y}}_k} - \underbrace{\begin{bmatrix} 
        \epsilon^{\frac{N-k}{2}}M_x \\ 
        M_T 
    \end{bmatrix}}_{C_{k}} \bm{\xi}_k , \enskip k =1,2, \ldots,N,
    \end{array}
\end{equation}
where the constraint function \(\bm{G}(\bm{\xi}_k, \bm{\eta}_k)\colon \mathbb{R}^{n_x+2} \times \mathbb{R}^{n_u+1} \to \mathbb{R}^{n_G} \) is defined as
\begin{equation}
\bm{G}(\bm{\xi}_k, \bm{\eta}_k) \coloneqq 
\begin{bmatrix}
\begin{bmatrix}
0_{n_x}^\top & 1 & 0
\end{bmatrix} \bm{\xi}_k - \gamma \\
\bm{\eta}_k - \overline{\bm{u}}_{\max} \\
- \bm{\eta}_k + \overline{\bm{u}}_{\min}
\end{bmatrix},
\end{equation}
and the operator \([\bm{z}]_+\) is defined element-wise as
\(
[\bm{z}]_+ \coloneqq \max\{\bm{z}, 0_{n_G}\}
\), where \(n_G = 2n_u + 3\). Here \(\nu \in \mathbb{R}_{>0}\) is a weighting parameter. By choosing an appropriate value for \(\nu\), the penalized term \(\norm{[\bm{G}(\bm{\xi}_k,\bm{\eta}_k)]_+ + \nu \ones_{n_G}}_2^2\) ensures that \( \bm{G}(\bm{\xi}_k,\bm{\eta}_k) \leq 0_{n_G}\)~\cite[Thm.~17.4]{nocedal1999numerical}. Positive-definite matrices \(Q \in \mathbb{R}^{6m \times 6m}\) and \(R \in \mathbb{R}^{(n_u+1) \times (n_u+1)}\) set the relative weighting on state and input tracking, respectively. The specific values of \(Q\) and \(R\) are selected in Section~\ref{sec:numerical_experiments}.

The idea of constraint-aware particle filtering is to transform the optimization~\eqref{pf: MHE_before} into an equivalent Bayesian state estimation problem. In particular, by introducing auxiliary variables, we can write optimization~\eqref{pf: MHE_before} equivalently as
\begin{equation}\label{Ed: min time relaxed}
    \begin{array}{ll}
    \underset{\substack{\bm{\xi}_{1:N}, \bm{\eta}_{1:N} \\
    \bm{\vartheta}_{0:N-1}, \bm{e}_{1:N}, \bm{h}_{1:N}}}{\mbox{minimize}} & \sum_{k=1}^{N} \big( \norm{ \bm{e}_k }_Q^2 + \norm{\bm{h}_k}_2^2 + \norm{\bm{\vartheta}_{k-1}}_{R}^2 \big) \\
    \mbox{subject to} & 
    \begin{bmatrix}
    \bm{\xi}_1 \\ \bm{\eta}_1  
    \end{bmatrix} = \begin{bmatrix}
        \overline{\bm{x}}_0 \\
        \hat{\bm{\eta}}
    \end{bmatrix} + \begin{bmatrix}
        0_{n_x+2} \\ 
        \bm{\vartheta}_0
    \end{bmatrix}, \\
    & \begin{bmatrix}
    \bm{\xi}_{k+1} \\ \bm{\eta}_{k+1} 
    \end{bmatrix} = \underbrace{\begin{bmatrix}
         \bm{F}_k(\bm{\xi}_k, \bm{\eta}_k) \\ 
         \hat{\bm{\eta}}
    \end{bmatrix}}_{\bm{\phi}(\bm{\chi}_k)} + \begin{bmatrix}
        0_{n_x+2} \\ 
        \bm{\vartheta}_k
    \end{bmatrix},   \enskip k = 1, 2, \ldots, N-1,\\
    & \underbrace{\begin{bmatrix}
    \hat{\bm{y}}_{k} \\ - \nu \bm{1}_{n_{G}}
    \end{bmatrix}}_{\tilde{\bm{y}}_k} = \underbrace{\begin{bmatrix}
         C_k \bm{\xi}_k \\ 
         [\bm{G}(\bm{\xi}_k,\bm{\eta}_k)]_+
    \end{bmatrix}}_{\bm{\psi}(\bm{\chi}_k)} + \begin{bmatrix}
        \bm{e}_k \\ 
        \bm{h}_k
    \end{bmatrix}, \enskip k = 1, 2, \ldots, N,
    \end{array}
\end{equation}
where \(\hat{\bm{\eta}} = \begin{bmatrix}
    0_{n_u}^\top & t_{\min}
\end{bmatrix}^\top \in \mathbb{R}^{n_u+1}\), \( \bm{\chi}_k = \begin{bmatrix}
   \bm{\xi}_k^\top & \bm{\eta}_k^\top
\end{bmatrix}^\top \in \mathbb{R}^{n_x+n_u+3}\). This transformation introduces \(\{\bm{h}_k\}_{k=1}^N\)---which penalizes all inequality constraint violations---and \(\{\bm{\vartheta}_k\}_{k=0}^{N-1}\), representing the original input tracking shifted one step backward to ensure temporal consistency in the update of the combined state \(\bm{\chi}_{k+1}\).


One can show that optimization~\eqref{Ed: min time relaxed} is equivalent to the following estimation problem~\cite[Thm. 1]{askari2023model}:
\begin{equation}\label{eqn: max likelihood}
    \begin{array}{ll}
    \underset{\bm{\chi_{1:N}}}{\mbox{maximize}}  & \log  \big( p(\bm{\chi}_{1:N} | \bm{y}_{1:N}=\tilde{\bm{y}}_{1:N}) \big), \\
    \end{array}
\end{equation}
where \(\bm{\chi_{1:N}}\) and \(\bm{y}_{1:N}\) satisfy the following dynamics equations: 
\begin{subequations}
    \begin{align}
 \bm{\chi}_{k+1}&=\bm{\phi}(\bm{\chi}_k)+\begin{bmatrix}
        0_{n_x+2}\\
        \bm{\vartheta}_k
    \end{bmatrix},\enskip \bm{\vartheta}_k \sim \mathcal{N}(0_{n_u+1}, R^{-1}),\\
    \bm{y}_k & = \bm{\psi} (\bm{\chi_k})+\begin{bmatrix}
        \bm{e}_k\\
        \bm{h}_k
    \end{bmatrix},\enskip \begin{bmatrix}
        \bm{e}_k\\
        \bm{h}_k
    \end{bmatrix}\sim \mathcal{N} (0_{6m+n_{G}}, \blkdiag(Q^{-1}, I_{n_{G}})).
    \end{align}
\end{subequations}
We define the process noise covariance \(E\) and the measurement noise covariance \(F\) as
\begin{equation}\label{Covariance matrix E and F}
    E \coloneqq \blkdiag(\bm{0}_{n_x+2}, R^{-1}), \enskip F \coloneqq \blkdiag(Q^{-1}, I_{n_G}).
\end{equation}


\subsection{Particle Filtering and Particle Selection}
The idea of constraint-aware particle filtering is to approximate the maximum likelihood distribution of \(\{\bm{\chi}_k\}_{k=1}^N\)~\eqref{eqn: max likelihood} using a finite number of particles~\cite{askari2023model}. Intuitively, the particles with the highest likelihood of satisfying both the system dynamics and the reference tracking in~\eqref{Ed: min time relaxed} can serve as approximate solutions to the optimization problem in~\eqref{pf: MHE_before}. We summarize the constraint-aware particle filtering method in Algorithm~\ref{alg: filter}, which uses the unscented transform---which we summarize in Algorithm~\ref{alg: unscented}---to propagate the mean and covariance of the system~\cite{fang2018nonlinear}.

\begin{algorithm}
\caption{Unscented Transform~\cite{fang2018nonlinear}}
\label{alg: unscented}
\begin{algorithmic}[1]
\Require Mean \(x\in\mathbb{R}^n\), variance \( A_1\in\mathbb{S}_{\succ 0}^{n}, A_2\in\mathbb{S}_{\succeq 0}^l\), function \(\varphi:\mathbb{R}^n\to\mathbb{R}^l\), \(\theta \approx 0.1\).
\State \(\lambda\gets (\theta^2-1)n\)
\State \(\overline{X}\gets x\ones_{2n+1}^\top\)
\State \(X\gets \overline{X}+\sqrt{n+\lambda}\begin{bmatrix}
0_n & \sqrt{A_1} & -\sqrt{A_1}
\end{bmatrix}\)
\State \(a\gets\frac{1}{n+\lambda}\begin{bmatrix}
  \lambda & \frac{1}{2}\ones_{2n}^\top  
\end{bmatrix}^\top\)
\State \(b \gets \frac{1}{n+\lambda}\begin{bmatrix}
    \lambda +(n+\lambda)(3-\theta^2) & \frac{1}{2}\ones^\top_{2n}
\end{bmatrix}^\top\)
\State \(Y\gets\varphi(X)\) (\(\varphi\) applied column-wise)
\State \(y\gets Ya\)
\State\(\overline{Y}\gets y\ones_{2n+1}^\top\)
\State \(B_1 \gets (Y-\overline{Y})\diag(b)(Y-\overline{Y})^\top+A_2\)
\State \(B_2 \gets (X-\overline{X})\diag(b) (Y-\overline{Y})^\top\)
\Ensure \(y, B_1, B_2\)
\end{algorithmic}
\end{algorithm}

\begin{algorithm}
\caption{Constraint-Aware Parcile Filtering}
\label{alg: filter}
\begin{algorithmic}[1]
\Require Number of particles \(n_p\). Initial state \(\bm{\chi}_1^l\in\mathbb{R}^{n_x+n_u+3}\), initial state covariance \(\Sigma_1^l\in\mathbb{S}_{\succeq 0}^{n_x+n_u+3}\) for all \(l=1, 2, \ldots, n_p\). Process noise covariance \(E\in\mathbb{S}_{\succeq 0}^{n_x+n_u+3}\), measurement noise covariance \(F\in\mathbb{S}_{\succeq 0}^{6m+n_G}\). Desired output \(\tilde{\bm{y}}_k \in\mathbb{R}^{6m+n_G}\) for all \(k= 1, 2, \ldots, N\). Sampling parameter \(\alpha\in(0, 1)\) and \(\kappa\in(1, n_p)\).
\State Let \(\omega_k^l=\frac{1}{n_p}\) for all \( 1\leq l\leq n_p\) and \(1 \leq k \leq N.\)
\For{\(k=1, 2, \ldots, N-1\)}
\For{\(l=1, 2, \ldots, n_p\)}
\State \((\mu, M, N)\gets\text{Alg.~\ref{alg: unscented}}(\bm{\chi}_k^l, \Sigma_k^l, E, \bm{\phi})\) \Comment{UT for dynamics.}
\State \((\zeta, U, V)\gets\text{Alg.~\ref{alg: unscented}}(\mu, M, F, \bm{\psi})\) \Comment{UT for output.}
\State \(K\gets V (U)^{-1}\) \Comment{Compute Kalman Gain.}
\State \(\Sigma_{k+1}^l\gets M-KU K^\top\) \Comment{Update state covariance.}
\State \(z\sim \mathcal{N}(0_{n_x+n_u+3}, \alpha I_{n_x+n_u+3})\) \Comment{Generate  a random seed.}
\State \(\bm{\chi}_{k+1}^l\gets \mu+K(\tilde{\bm{y}}_{k+1}-\zeta)+\sqrt{\Sigma_{k+1}^l}z\) \Comment{Update state mean with random bias.}
\State \(\tilde{\omega}_{k+1}^l\gets\frac{\omega_k^l}{\sqrt{\det U}} \exp\left(-\frac{1}{2}\norm{\tilde{\bm{y}}_{k+1}-\zeta}_{(U)^{-1}}^2\right)\) \Comment{Evaluate particle likelihood.}
\EndFor
\State \(\omega_{k+1}^l \gets \tilde{\omega}_{k+1}^l/\left(\sum_{j=1}^{n_p}\tilde{\omega}_{k+1}^j\right)\) for all \(1\leq l \leq n_p.\)\Comment{Normalization.}
\If{\(\kappa\sum_{j=1}^{n_p}(\omega_{k+1}^j)^2\geq 1\) } \Comment{If needed, re-sample.}
\State Let \( \{\overline{\bm{\chi}}_s^l\}_{s=1}^{k+1} \gets \{\bm{\chi}_s^l\}_{s = 1}^{k+1}\) and \(\overline{\Sigma}_{k+1}^l\gets\Sigma_{k+1}^l\) for all \(1 \leq l \leq n_p\).
\For{\(l =1, 2, \ldots, n_p\)}
\State Sample a random integer \(j\) such that \(\mathds{P}(j=i)=\omega_{k+1}^i\) for all \(1\leq i \leq n_p\).
\State Let \(\{\bm{\chi}_s^l\}_{s=1}^{k+1}\gets\{\overline{\bm{\chi}}_s^j\}_{s=1}^{k+1}\), \(\Sigma_{k+1}^l\gets \overline{\Sigma}_{k+1}^j\) and  \(\omega_{k+1}^l\gets \frac{1}{m}\).
\EndFor 
\EndIf
\EndFor

\Ensure \(\{\{\bm{\chi}_k^l\}_{k=1}^N\}_{l=1}^{n_p}\). 
\end{algorithmic}
\end{algorithm}


Notice that the output of Algorithm~\ref{alg: filter} consists of \(n_p\) particle trajectories, where each trajectory \( \{\bm{\chi}_k^l\}_{k=1}^N\) represents the estimated states of particle \(l\) over \(N\) discretization steps. To initialize the iterative solution methods, such as the proposed SCP approach in Section~\ref{sec: proposed SCP}, we select the particle as the warm-starting trajectory based on the lowest objective function value and minimal state dynamics violation. We let
\begin{equation}
\begin{aligned}
    \phi^l = &\sum_{k=1}^{N}\left( \norm{\hat{\bm{y}}_k - C_k \bm{\xi}_k^l}_Q^2 + \norm{[\bm{G}(\bm{\xi}_k^l,\bm{\eta}_k^l)]_+ + \nu \ones_{n_G}}_2^2 + \norm{\bm{\eta}_k^l - \hat{\bm{\eta}}}_R^2\right)\\
    & + \sum_{k=1}^{N-1} \left( \norm{\bm{F}_k(\bm{\xi}_k^l,\bm{\eta}_k^l)- \bm{\xi}_{k+1}^l }_1
    \right),
\end{aligned}
\end{equation}
for all \(1\leq l \leq n_p\), where \(\bm{\xi}_k^l\) and \(\bm{\eta}_k^l\) are derived from \(\bm{\chi}_k^l\), with \(\bm{\chi}_k^l = \begin{bmatrix}
    \left(\bm{\xi}_k^l\right)^\top & \left(\bm{\eta}_k^l\right)^\top
\end{bmatrix}^\top\).
We define the 
\begin{equation}\label{eqn: good_particle}
    l^\star \coloneqq \argmin_{1 \leq l \leq n_p} \enskip \phi_l,
\end{equation}
and warm-start the Algorithm~\ref{alg: prox-linear} by letting 
\begin{equation}
    \{\hat{\bm{\xi}}_k^{(0)}\}_{k=1}^N = \{\bm{\xi}_k^{l^\star}\}_{k=1}^N,\enskip \{\hat{\bm{\eta}}_k^{(0)}\}_{k=1}^{N-1} =\{\bm{\eta}_k^{l^\star}\}_{k=1}^{N-1}.
\end{equation}

We base the warm-starting strategy above on previous work~\cite{yuan2024filtering}. However, the previous work selects the warm-starting trajectory from a set of weighted average of trajectories with the lowest objective function values and minimal state dynamics violations. To obtain these weighted average trajectories, the previous work relies on \emph{hierarchical clustering}~\cite[Sec.~8.2]{nielsen2016hierarchical} to identify distinct clusters in the sampled trajectories and then computes a weighted average trajectory within each cluster. In contrast, here we propose to select the warm-starting trajectory via~\eqref{eqn: good_particle} without performing \emph{hierarchical clustering}. Later we will show that this new selection still achieve good convergence performance in Section~\ref{sec:numerical_experiments}.





\section{Numerical Experiments}\label{sec:numerical_experiments}

We demonstrate the SCP approach proposed in Section~\ref{sec: proposed SCP} and the filtering-based warm-starting strategy proposed in Section~\ref{sec: proposed warm-starting} on multiagent quadrotor trajectory optimization problems with collision avoidance constraints. We also compare its convergence performance with \texttt{IPOPT}, an off-the-shelf interior-point method solver~\cite{wachter2006implementation}.



\subsection{Problem Setup}

We consider a multiagent quadrotor trajectory optimization problem in a shared space containing cylindrical obstacles. Below, we list the parameters for optimization~\eqref{ocp: min energy and time} as well as the parameters for Algorithm~\ref{alg: prox-linear} and Algorithm~\ref{alg: filter}.

\subsubsection{Problem Parameters}

We set the parameters in optimization~\eqref{ocp: min energy and time} as follows. We consider the cases \(m = 2\), \(m = 4\), and \(m = 6\), all with \(n_o = 2\). We list the quadrotor mass in Section~\ref{subsection: Dynamics}, the normalized weighting parameters in optimization~\eqref{ocp: min energy and time}, the path constraint-related parameters in Section~\ref{subsection: Path Constraints}, and the final time constraint parameters in Section~\ref{subsection: Boundary Constraints} in Table~\ref{tab:Problem parameters} for all three cases. In each case, we set the initial and final velocities in \(\hat{\bm{x}}_0^i\) and \(\hat{\bm{x}}_f^i\) to \(0_3\), and the initial and final thrusts to \(\hat{\bm{T}}\), as defined in Section~\ref{subsection: Dynamics}, for all agents. We select the initial and final positions of all agents from the vertices and edge midpoints of a bounding box defined by two diagonally opposite corners with coordinates \( \begin{bmatrix}
     2 & 2 &2
\end{bmatrix}^\top\,\mathrm{m}\) and \(\begin{bmatrix}
    14 & 14 & 14
\end{bmatrix}^\top\,\mathrm{m}\), ensuring that the agents start and end at spatially well-separated locations. Notice that the minimum and maximum time bounds in Table~\ref{tab:Problem parameters} are based on the spatial distribution of these positions and the maximum velocity.


\begin{table}[hbt!]
\caption{\label{tab:Problem parameters} Parameter values in optimization~\eqref{ocp: min energy and time}}
\centering
\begin{tabular}{lclc}
\toprule
Parameter & Value & Parameter & Value \\
\midrule
$\bm{c}_1,\,\mathrm{m}$   & $\begin{bmatrix}  5 & 8 \end{bmatrix}^\top$   &    $\rho_1,\,\mathrm{m}$   & 2 \\
$\bm{c}_2,\,\mathrm{m}$   & $\begin{bmatrix}  9 & 5 \end{bmatrix}^\top$   &    $\rho_2,\,\mathrm{m}$   & 1.5 \\
$d,\,\mathrm{m}$   & 1   & $m_q,\,\mathrm{kg}$   & 0.35       \\
$\bm{r}_{\max},\,\mathrm{m}$  & $\begin{bmatrix}15 & 15 & 15\end{bmatrix}^\top $  & $\bm{r}_{\min},\,\mathrm{m}$  & $\begin{bmatrix}0 & 0 & 0\end{bmatrix}^\top $   \\
$T_{\max},\,\mathrm{N}$    & 5.00    & $T_{\min},\,\mathrm{N}$    & 2.00        \\
$\bm{u}_{\max},\,\mathrm{N/s}$  & $\begin{bmatrix}2 & 2 & 2\end{bmatrix}^\top $   & $\bm{u}_{\min},\,\mathrm{N/s}$  & $\begin{bmatrix}-2 & -2 & -2\end{bmatrix}^\top $  \\
$v_{\max},\,\mathrm{m/s}$  & 3.00   & $\theta_{\max},\,\mathrm{rad}$  & $\frac{\pi}{4}$       \\
$t_{\min},\,\mathrm{s}$  & 7   & $t_{\max},\,\mathrm{s}$  & 28     \\
$\tilde{\alpha}_1,\,\mathrm{unitless}$  & 0.1   & $\tilde{\alpha}_2,\,\mathrm{unitless}$  & 0.8     \\
$\tilde{\alpha}_3,\,\mathrm{unitless}$  & 0.1   &    \\
\bottomrule
\end{tabular}
\end{table}

\subsubsection{Algorithm parameters}

\paragraph{Prox-linear method:} Since we use adaptive-time discretization and account for constraints violation in continuous time, a relatively small number of discretization points is sufficient. Specifically, we pick \(N = 8\), which includes both the initial and final time steps. We set the parameter \(\gamma\) in optimization~\eqref{opt: min_time_relaxed_prox} to \(10^{-6}\), and set the parameters \(\beta = 20\) and \(\rho = 0.1\) in Algorithm~\ref{alg: prox-linear}. We terminate the Algorithm~\ref{alg: prox-linear} after a fixed runtime.

\paragraph{Particle filtering:} We choose the weighting matrices \(Q = I_{6m}\) and \(R = \blkdiag\left(\frac{\alpha_2}{\alpha_3} I_{n_u}, \frac{\alpha_1}{\alpha_3 t_{\max}}\right)\), and then obtain the process noise covariance \(E\) and measurement noise covariance \(F\) from~\eqref{Covariance matrix E and F}. We set \(n_p = 30\), where \(n_p\) denotes the number of particles. The initial state and initial state covariance are set to \(\bm{\chi}_1^l = \begin{bmatrix} \overline{\bm{x}}_0^\top & \hat{\bm{\eta}}^\top \end{bmatrix}^\top\) and \(\Sigma_1^l = 10^{-2} \cdot I_{n_x + n_u + 3}\), respectively, for all \(1 \leq l \leq n_p\). We set the scalar parameters \(\epsilon\) and \(\nu\) in optimization~\eqref{pf: MHE_before} to 0.5 and 1, respectively. We set the sampling parameters \(\alpha\) and \(\kappa\) in Algorithm~\ref{alg: filter} to \(5 \times 10^{-3}\) and 9, respectively.

\subsection{Experimental Setup}

We first compare the convergence performance of the proposed SCP approach in Section~\ref{sec: proposed SCP} with that of the state‑of‑the‑art nonlinear programming solver \texttt{IPOPT}~\cite{wachter2006implementation} (using both multiple‑shooting and single‑shooting formulations), and run all methods under random initialization. We then apply the warm‑starting strategy proposed in Section~\ref{sec: proposed warm-starting} to all methods to evaluate its impact on their performance. We implement the proposed SCP approach and all benchmark methods in \texttt{MATLAB}, and run simulations on the Minnesota Supercomputing Institute’s resources\footnote{\url{https://www.msi.umn.edu/}}. For the proposed SCP approach, we use the \texttt{OSQP}~\cite{osqp} solver to solve the convex subproblem~\eqref{opt: prox-linear iter} in Algorithm~\ref{alg: prox-linear}. As for the benchmark method, we use \texttt{CasADi}~\cite{andersson2019casadi} as the symbolic framework to model system dynamics and solve the optimization~\eqref{opt: min_time_relaxed_v1} with \texttt{IPOPT} via \texttt{CasADi}’s built‑in interface. We access both \texttt{OSQP} and \texttt{IPOPT}, which are implemented in \texttt{C}/\texttt{C++}, through their respective \texttt{MATLAB} interfaces. We also employ \texttt{MATLAB} \texttt{codegen} to convert the warm-start generation code to \texttt{C}/\texttt{C++}.


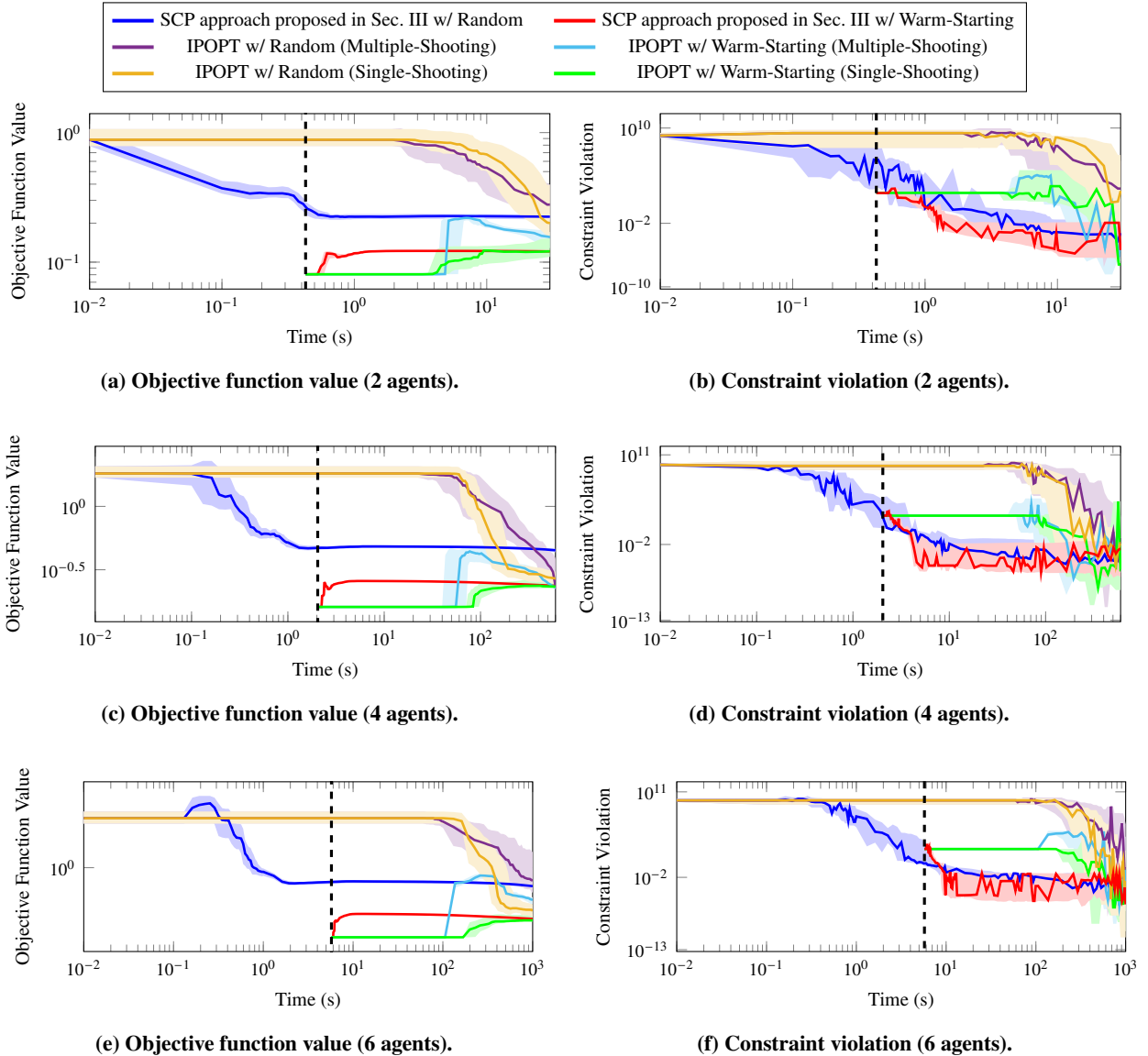
\begin{figure}[!htp]
\centering

\begin{tikzpicture}[baseline]
\pgfplotslegendfromname{sharedlegend}
\end{tikzpicture}

\begin{subfigure}[t]{0.5\linewidth}
\centering
\pgfplotsset{every tick label/.append style={font=\scriptsize}}
\begin{tikzpicture}
\begin{loglogaxis}[
    xlabel = \footnotesize Time (s), ylabel=\footnotesize Objective Function Value,
    width= \linewidth,
    height = 0.5\linewidth,
    xmin = 1e-2, xmax = 30,
    legend to name = sharedlegend,          
    legend style={
        legend columns=2,
        /tikz/column 2/.style={
            column sep=10pt},
        font=\footnotesize
    },
    legend image post style={line width=1.5pt} 
]


\definecolor{proxrandom_color}{rgb}{0,0,1}
\addplot [mark=none, color=proxrandom_color, line width=1pt] table[x=averageTimes, y=objMedianVals] {2agent_proxrandom.dat};
\addlegendentry{SCP approach proposed in Sec.~\ref{sec: proposed SCP} w/ Random}

\definecolor{proxws_color}{rgb}{1,0,0}
\addplot [mark=none, color=proxws_color,line width=1pt] table[x=averageTimes,y=objMedianVals] {2agent_proxws.dat};
\addlegendentry{SCP approach proposed in Sec.~\ref{sec: proposed SCP} w/ Warm-Starting}

\definecolor{ipoptran_MS_color}{rgb}{0.4940,0.1840,0.5560}
\addplot [mark=none, color=ipoptran_MS_color,line width=1pt] table[x=averageTimes,y=objMedianVals] {2agent_ipoptran_MS.dat};
\addlegendentry{IPOPT w/ Random (Multiple-Shooting)}

\definecolor{ipoptws_MS_color}{rgb}{0.3010,0.7450,0.9330}
\addplot [mark=none, color=ipoptws_MS_color,line width=1pt] table[x=averageTimes,y=objMedianVals] {2agent_ipoptws_MS.dat};
\addlegendentry{IPOPT w/ Warm-Starting (Multiple-Shooting)}

\definecolor{ipoptran_SS_color}{rgb}{0.9290,0.6940,0.1250}
\addplot [mark=none, color=ipoptran_SS_color,line width=1pt] table[x=averageTimes,y=objMedianVals] {2agent_ipoptran_SS.dat};
\addlegendentry{IPOPT w/ Random (Single-Shooting)}

\definecolor{ipoptws_SS_color}{rgb}{0,1,0}
\addplot [mark=none, color=ipoptws_SS_color,line width=1pt] table[x=averageTimes,y=objMedianVals] {2agent_ipoptws_SS.dat};
\addlegendentry{IPOPT w/ Warm-Starting (Single-Shooting)}


\addplot [name path=upper, draw=none] table[x=averageTimes, y=objUpperMedianVals] {2agent_proxrandom.dat};
\addplot [name path=lower, draw=none] table[x=averageTimes, y=objLowerMedianVals] {2agent_proxrandom.dat};
\addplot [fill=proxrandom_color!20] fill between[of=upper and lower];


\addplot [name path=upper,draw=none] table[x=averageTimes,y = objUpperMedianVals] {2agent_proxws.dat};
\addplot [name path=lower,draw=none] table[x=averageTimes,y = objLowerMedianVals] {2agent_proxws.dat};
\addplot [fill=proxws_color!20] fill between[of=upper and lower];


\addplot [name path=upper,draw=none] table[x=averageTimes,y = objUpperMedianVals] {2agent_ipoptws_MS.dat};
\addplot [name path=lower,draw=none] table[x=averageTimes,y = objLowerMedianVals] {2agent_ipoptws_MS.dat};
\addplot [fill=ipoptws_MS_color!20] fill between[of=upper and lower];


\addplot [name path=upper,draw=none] table[x=averageTimes,y = objUpperMedianVals] {2agent_ipoptran_MS.dat};
\addplot [name path=lower,draw=none] table[x=averageTimes,y = objLowerMedianVals] {2agent_ipoptran_MS.dat};
\addplot [fill=ipoptran_MS_color!20] fill between[of=upper and lower];

\addplot [name path=upper,draw=none] table[x=averageTimes,y = objUpperMedianVals] {2agent_ipoptws_SS.dat};
\addplot [name path=lower,draw=none] table[x=averageTimes,y = objLowerMedianVals] {2agent_ipoptws_SS.dat};
\addplot [fill=ipoptws_SS_color!20] fill between[of=upper and lower];

\addplot [name path=upper,draw=none] table[x=averageTimes,y = objUpperMedianVals] {2agent_ipoptran_SS.dat};
\addplot [name path=lower,draw=none] table[x=averageTimes,y = objLowerMedianVals] {2agent_ipoptran_SS.dat};
\addplot [fill=ipoptran_SS_color!20] fill between[of=upper and lower];

\draw[dashed, very thick, color=black] (axis cs:0.428346375000000, 1e-3) -- (axis cs:0.428346375000000, 1e3);

\end{loglogaxis}
\end{tikzpicture}
\caption{Objective function value (2 agents).}    
\end{subfigure}\hfill
\begin{subfigure}[t]{0.5\linewidth}
\centering
\pgfplotsset{every tick label/.append style={font=\scriptsize}}
\begin{tikzpicture}
\begin{loglogaxis}[
    xlabel = \footnotesize Time (s), ylabel=\footnotesize Constraint Violation,
    width= \linewidth,
    height = 0.5\linewidth,
    xmin = 1e-2, xmax = 30,
    ytick={1e-10, 1e-2, 1e10}
]

\definecolor{proxrandom_color}{rgb}{0,0,1}
\addplot [mark=none, color=proxrandom_color,line width=1pt] table[x=averageTimes,y=conMedianVals] {2agent_proxrandom.dat};

\addplot [name path=upper,draw=none] table[x=averageTimes,y = conUpperMedianVals] {2agent_proxrandom.dat};
\addplot [name path=lower,draw=none] table[x=averageTimes,y = conLowerMedianVals] {2agent_proxrandom.dat};
\addplot [fill=proxrandom_color!20] fill between[of=upper and lower];


\definecolor{ipoptws_MS_color}{rgb}{0.301, 0.745, 0.933}
\addplot [mark=none, color=ipoptws_MS_color,line width=1pt] table[x=averageTimes,y=conMedianVals] {2agent_ipoptws_MS.dat};

\addplot [name path=upper,draw=none] table[x=averageTimes,y = conUpperMedianVals] {2agent_ipoptws_MS.dat};
\addplot [name path=lower,draw=none] table[x=averageTimes,y = conLowerMedianVals] {2agent_ipoptws_MS.dat};
\addplot [fill=ipoptws_MS_color!20] fill between[of=upper and lower];

\definecolor{ipoptran_MS_color}{rgb}{0.4940,0.1840,0.5560}
\addplot [mark=none, color=ipoptran_MS_color,line width=1pt] table[x=averageTimes,y=conMedianVals] {2agent_ipoptran_MS.dat};

\addplot [name path=upper,draw=none] table[x=averageTimes,y = conUpperMedianVals] {2agent_ipoptran_MS.dat};
\addplot [name path=lower,draw=none] table[x=averageTimes,y = conLowerMedianVals] {2agent_ipoptran_MS.dat};
\addplot [fill=ipoptran_MS_color!20] fill between[of=upper and lower];

\definecolor{ipoptws_SS_color}{rgb}{0, 1, 0}
\addplot [mark=none, color=ipoptws_SS_color,line width=1pt] table[x=averageTimes,y=conMedianVals] {2agent_ipoptws_SS.dat};

\addplot [name path=upper,draw=none] table[x=averageTimes,y = conUpperMedianVals] {2agent_ipoptws_SS.dat};
\addplot [name path=lower,draw=none] table[x=averageTimes,y = conLowerMedianVals] {2agent_ipoptws_SS.dat};
\addplot [fill=ipoptws_SS_color!20] fill between[of=upper and lower];

\definecolor{ipoptran_SS_color}{rgb}{0.9290,0.6940,0.1250}
\addplot [mark=none, color=ipoptran_SS_color,line width=1pt] table[x=averageTimes,y=conMedianVals] {2agent_ipoptran_SS.dat};
\addplot [name path=upper,draw=none] table[x=averageTimes,y = conUpperMedianVals] {2agent_ipoptran_SS.dat};
\addplot [name path=lower,draw=none] table[x=averageTimes,y = conLowerMedianVals] {2agent_ipoptran_SS.dat};
\addplot [fill=ipoptran_SS_color!20] fill between[of=upper and lower];

\definecolor{proxws_color}{rgb}{1,0,0}
\addplot [mark=none, color=proxws_color,line width=1pt] table[x=averageTimes,y=conMedianVals] {2agent_proxws.dat};

\addplot [name path=upper,draw=none] table[x=averageTimes,y = conUpperMedianVals] {2agent_proxws.dat};
\addplot [name path=lower,draw=none] table[x=averageTimes,y = conLowerMedianVals] {2agent_proxws.dat};
\addplot [fill=proxws_color!20] fill between[of=upper and lower];

\draw[dashed, very thick, color=black] (axis cs:0.428346375000000, 1e-11) -- (axis cs:0.428346375000000, 1e12);

\end{loglogaxis}
\end{tikzpicture}
\caption{Constraint violation (2 agents).}
\end{subfigure}

\bigskip

\begin{subfigure}[t]{0.5\linewidth}
\centering
\pgfplotsset{every tick label/.append style={font=\scriptsize}}
\begin{tikzpicture}
\begin{loglogaxis}[
    xlabel = \footnotesize Time (s), ylabel=\footnotesize Objective Function Value,
    width= \linewidth,
    height = 0.5\linewidth,
    xmin = 1e-2, xmax = 600,
]


\definecolor{proxrandom_color}{rgb}{0,0,1}
\addplot [mark=none, color=proxrandom_color, line width=1pt] table[x=averageTimes, y=objMedianVals] {4agent_proxrandom.dat};

\definecolor{proxws_color}{rgb}{1,0,0}
\addplot [mark=none, color=proxws_color,line width=1pt] table[x=averageTimes,y=objMedianVals] {4agent_proxws.dat};

\definecolor{ipoptran_MS_color}{rgb}{0.4940,0.1840,0.5560}
\addplot [mark=none, color=ipoptran_MS_color,line width=1pt] table[x=averageTimes,y=objMedianVals] {4agent_ipoptran_MS.dat};

\definecolor{ipoptws_MS_color}{rgb}{0.3010,0.7450,0.9330}
\addplot [mark=none, color=ipoptws_MS_color,line width=1pt] table[x=averageTimes,y=objMedianVals] {4agent_ipoptws_MS.dat};

\definecolor{ipoptran_SS_color}{rgb}{0.9290,0.6940,0.1250}
\addplot [mark=none, color=ipoptran_SS_color,line width=1pt] table[x=averageTimes,y=objMedianVals] {4agent_ipoptran_SS.dat};

\definecolor{ipoptws_SS_color}{rgb}{0,1,0}
\addplot [mark=none, color=ipoptws_SS_color,line width=1pt] table[x=averageTimes,y=objMedianVals] {4agent_ipoptws_SS.dat};


\addplot [name path=upper, draw=none] table[x=averageTimes, y=objUpperMedianVals] {4agent_proxrandom.dat};
\addplot [name path=lower, draw=none] table[x=averageTimes, y=objLowerMedianVals] {4agent_proxrandom.dat};
\addplot [fill=proxrandom_color!20] fill between[of=upper and lower];


\addplot [name path=upper,draw=none] table[x=averageTimes,y = objUpperMedianVals] {4agent_proxws.dat};
\addplot [name path=lower,draw=none] table[x=averageTimes,y = objLowerMedianVals] {4agent_proxws.dat};
\addplot [fill=proxws_color!20] fill between[of=upper and lower];


\addplot [name path=upper,draw=none] table[x=averageTimes,y = objUpperMedianVals] {4agent_ipoptws_MS.dat};
\addplot [name path=lower,draw=none] table[x=averageTimes,y = objLowerMedianVals] {4agent_ipoptws_MS.dat};
\addplot [fill=ipoptws_MS_color!20] fill between[of=upper and lower];


\addplot [name path=upper,draw=none] table[x=averageTimes,y = objUpperMedianVals] {4agent_ipoptran_MS.dat};
\addplot [name path=lower,draw=none] table[x=averageTimes,y = objLowerMedianVals] {4agent_ipoptran_MS.dat};
\addplot [fill=ipoptran_MS_color!20] fill between[of=upper and lower];

\addplot [name path=upper,draw=none] table[x=averageTimes,y = objUpperMedianVals] {4agent_ipoptws_SS.dat};
\addplot [name path=lower,draw=none] table[x=averageTimes,y = objLowerMedianVals] {4agent_ipoptws_SS.dat};
\addplot [fill=ipoptws_SS_color!20] fill between[of=upper and lower];

\addplot [name path=upper,draw=none] table[x=averageTimes,y = objUpperMedianVals] {4agent_ipoptran_SS.dat};
\addplot [name path=lower,draw=none] table[x=averageTimes,y = objLowerMedianVals] {4agent_ipoptran_SS.dat};
\addplot [fill=ipoptran_SS_color!20] fill between[of=upper and lower];

\draw[dashed, very thick, color=black] (axis cs:2.04093973300000, 1e-3) -- (axis cs:2.04093973300000, 1e3);

\end{loglogaxis}
\end{tikzpicture}
\caption{Objective function value (4 agents).}    
\end{subfigure}\hfill
\begin{subfigure}[t]{0.5\linewidth}
\centering
\pgfplotsset{every tick label/.append style={font=\scriptsize}}
\begin{tikzpicture}
\begin{loglogaxis}[
    xlabel = \footnotesize Time (s), ylabel=\footnotesize Constraint Violation,
    width= \linewidth,
    height = 0.5\linewidth,
    xmin = 1e-2, xmax = 600,
    ytick = {1e-13, 1e-2, 1e11}
]

\definecolor{proxrandom_color}{rgb}{0,0,1}
\addplot [mark=none, color=proxrandom_color,line width=1pt] table[x=averageTimes,y=conMedianVals] {4agent_proxrandom.dat};

\addplot [name path=upper,draw=none] table[x=averageTimes,y = conUpperMedianVals] {4agent_proxrandom.dat};
\addplot [name path=lower,draw=none] table[x=averageTimes,y = conLowerMedianVals] {4agent_proxrandom.dat};
\addplot [fill=proxrandom_color!20] fill between[of=upper and lower];


\definecolor{ipoptws_MS_color}{rgb}{0.301, 0.745, 0.933}
\addplot [mark=none, color=ipoptws_MS_color,line width=1pt] table[x=averageTimes,y=conMedianVals] {4agent_ipoptws_MS.dat};

\addplot [name path=upper,draw=none] table[x=averageTimes,y = conUpperMedianVals] {4agent_ipoptws_MS.dat};
\addplot [name path=lower,draw=none] table[x=averageTimes,y = conLowerMedianVals] {4agent_ipoptws_MS.dat};
\addplot [fill=ipoptws_MS_color!20] fill between[of=upper and lower];

\definecolor{ipoptran_MS_color}{rgb}{0.4940,0.1840,0.5560}
\addplot [mark=none, color=ipoptran_MS_color,line width=1pt] table[x=averageTimes,y=conMedianVals] {4agent_ipoptran_MS.dat};

\addplot [name path=upper,draw=none] table[x=averageTimes,y = conUpperMedianVals] {4agent_ipoptran_MS.dat};
\addplot [name path=lower,draw=none] table[x=averageTimes,y = conLowerMedianVals] {4agent_ipoptran_MS.dat};
\addplot [fill=ipoptran_MS_color!20] fill between[of=upper and lower];

\definecolor{ipoptws_SS_color}{rgb}{0, 1, 0}
\addplot [mark=none, color=ipoptws_SS_color,line width=1pt] table[x=averageTimes,y=conMedianVals] {4agent_ipoptws_SS.dat};

\addplot [name path=upper,draw=none] table[x=averageTimes,y = conUpperMedianVals] {4agent_ipoptws_SS.dat};
\addplot [name path=lower,draw=none] table[x=averageTimes,y = conLowerMedianVals] {4agent_ipoptws_SS.dat};
\addplot [fill=ipoptws_SS_color!20] fill between[of=upper and lower];

\definecolor{ipoptran_SS_color}{rgb}{0.9290,0.6940,0.1250}
\addplot [mark=none, color=ipoptran_SS_color,line width=1pt] table[x=averageTimes,y=conMedianVals] {4agent_ipoptran_SS.dat};
\addplot [name path=upper,draw=none] table[x=averageTimes,y = conUpperMedianVals] {4agent_ipoptran_SS.dat};
\addplot [name path=lower,draw=none] table[x=averageTimes,y = conLowerMedianVals] {4agent_ipoptran_SS.dat};
\addplot [fill=ipoptran_SS_color!20] fill between[of=upper and lower];

\definecolor{proxws_color}{rgb}{1,0,0}
\addplot [mark=none, color=proxws_color,line width=1pt] table[x=averageTimes,y=conMedianVals] {4agent_proxws.dat};

\addplot [name path=upper,draw=none] table[x=averageTimes,y = conUpperMedianVals] {4agent_proxws.dat};
\addplot [name path=lower,draw=none] table[x=averageTimes,y = conLowerMedianVals] {4agent_proxws.dat};
\addplot [fill=proxws_color!20] fill between[of=upper and lower];

\draw[dashed, very thick, color=black] (axis cs:2.04093973300000, 1e-13) -- (axis cs:2.04093973300000, 1e12);

\end{loglogaxis}
\end{tikzpicture}
\caption{Constraint violation (4 agents).}
\end{subfigure}

\bigskip

\begin{subfigure}[t]{0.49\linewidth}
\centering
\raisebox{0pt}{
\pgfplotsset{every tick label/.append style={font=\scriptsize}}
\begin{tikzpicture}
\begin{loglogaxis}[
    xlabel = \footnotesize Time (s), ylabel=\footnotesize Objective Function Value,
    width= \linewidth,
    height = 0.5\linewidth,
    xmin = 1e-2, xmax = 1000
]


\definecolor{proxrandom_color}{rgb}{0,0,1}
\addplot [mark=none, color=proxrandom_color, line width=1pt] table[x=averageTimes, y=objMedianVals] {6agent_proxrandom.dat};

\definecolor{proxws_color}{rgb}{1,0,0}
\addplot [mark=none, color=proxws_color,line width=1pt] table[x=averageTimes,y=objMedianVals] {6agent_proxws.dat};

\definecolor{ipoptran_MS_color}{rgb}{0.4940,0.1840,0.5560}
\addplot [mark=none, color=ipoptran_MS_color,line width=1pt] table[x=averageTimes,y=objMedianVals] {6agent_ipoptran_MS.dat};

\definecolor{ipoptws_MS_color}{rgb}{0.3010,0.7450,0.9330}
\addplot [mark=none, color=ipoptws_MS_color,line width=1pt] table[x=averageTimes,y=objMedianVals] {6agent_ipoptws_MS.dat};

\definecolor{ipoptran_SS_color}{rgb}{0.9290,0.6940,0.1250}
\addplot [mark=none, color=ipoptran_SS_color,line width=1pt] table[x=averageTimes,y=objMedianVals] {6agent_ipoptran_SS.dat};

\definecolor{ipoptws_SS_color}{rgb}{0,1,0}
\addplot [mark=none, color=ipoptws_SS_color,line width=1pt] table[x=averageTimes,y=objMedianVals] {6agent_ipoptws_SS.dat};


\addplot [name path=upper, draw=none] table[x=averageTimes, y=objUpperMedianVals] {6agent_proxrandom.dat};
\addplot [name path=lower, draw=none] table[x=averageTimes, y=objLowerMedianVals] {6agent_proxrandom.dat};
\addplot [fill=proxrandom_color!20] fill between[of=upper and lower];


\addplot [name path=upper,draw=none] table[x=averageTimes,y = objUpperMedianVals] {6agent_proxws.dat};
\addplot [name path=lower,draw=none] table[x=averageTimes,y = objLowerMedianVals] {6agent_proxws.dat};
\addplot [fill=proxws_color!20] fill between[of=upper and lower];


\addplot [name path=upper,draw=none] table[x=averageTimes,y = objUpperMedianVals] {6agent_ipoptws_MS.dat};
\addplot [name path=lower,draw=none] table[x=averageTimes,y = objLowerMedianVals] {6agent_ipoptws_MS.dat};
\addplot [fill=ipoptws_MS_color!20] fill between[of=upper and lower];


\addplot [name path=upper,draw=none] table[x=averageTimes,y = objUpperMedianVals] {6agent_ipoptran_MS.dat};
\addplot [name path=lower,draw=none] table[x=averageTimes,y = objLowerMedianVals] {6agent_ipoptran_MS.dat};
\addplot [fill=ipoptran_MS_color!20] fill between[of=upper and lower];

\addplot [name path=upper,draw=none] table[x=averageTimes,y = objUpperMedianVals] {6agent_ipoptws_SS.dat};
\addplot [name path=lower,draw=none] table[x=averageTimes,y = objLowerMedianVals] {6agent_ipoptws_SS.dat};
\addplot [fill=ipoptws_SS_color!20] fill between[of=upper and lower];

\addplot [name path=upper,draw=none] table[x=averageTimes,y = objUpperMedianVals] {6agent_ipoptran_SS.dat};
\addplot [name path=lower,draw=none] table[x=averageTimes,y = objLowerMedianVals] {6agent_ipoptran_SS.dat};
\addplot [fill=ipoptran_SS_color!20] fill between[of=upper and lower];

\draw[dashed, very thick, color=black] (axis cs:5.73827373350000, 1e-3) -- (axis cs:5.73827373350000, 1e3);

\end{loglogaxis}
\end{tikzpicture}}
\caption{Objective function value (6 agents).}    
\end{subfigure}\hfill
\begin{subfigure}[t]{0.49\linewidth}
\centering
\raisebox{0pt}{
\pgfplotsset{every tick label/.append style={font=\scriptsize}}
\begin{tikzpicture}
\begin{loglogaxis}[
    xlabel = \footnotesize Time (s), ylabel=\footnotesize Constraint Violation,
    width= \linewidth,
    height = 0.5\linewidth,
    xmin = 1e-2, xmax = 1000,
    ytick = {1e-13, 1e-2, 1e11}
]

\definecolor{proxrandom_color}{rgb}{0,0,1}
\addplot [mark=none, color=proxrandom_color,line width=1pt] table[x=averageTimes,y=conMedianVals] {6agent_proxrandom.dat};

\addplot [name path=upper,draw=none] table[x=averageTimes,y = conUpperMedianVals] {6agent_proxrandom.dat};
\addplot [name path=lower,draw=none] table[x=averageTimes,y = conLowerMedianVals] {6agent_proxrandom.dat};
\addplot [fill=proxrandom_color!20] fill between[of=upper and lower];


\definecolor{ipoptws_MS_color}{rgb}{0.301, 0.745, 0.933}
\addplot [mark=none, color=ipoptws_MS_color,line width=1pt] table[x=averageTimes,y=conMedianVals] {6agent_ipoptws_MS.dat};

\addplot [name path=upper,draw=none] table[x=averageTimes,y = conUpperMedianVals] {6agent_ipoptws_MS.dat};
\addplot [name path=lower,draw=none] table[x=averageTimes,y = conLowerMedianVals] {6agent_ipoptws_MS.dat};
\addplot [fill=ipoptws_MS_color!20] fill between[of=upper and lower];

\definecolor{ipoptran_MS_color}{rgb}{0.4940,0.1840,0.5560}
\addplot [mark=none, color=ipoptran_MS_color,line width=1pt] table[x=averageTimes,y=conMedianVals] {6agent_ipoptran_MS.dat};

\addplot [name path=upper,draw=none] table[x=averageTimes,y = conUpperMedianVals] {6agent_ipoptran_MS.dat};
\addplot [name path=lower,draw=none] table[x=averageTimes,y = conLowerMedianVals] {6agent_ipoptran_MS.dat};
\addplot [fill=ipoptran_MS_color!20] fill between[of=upper and lower];

\definecolor{ipoptws_SS_color}{rgb}{0, 1, 0}
\addplot [mark=none, color=ipoptws_SS_color,line width=1pt] table[x=averageTimes,y=conMedianVals] {6agent_ipoptws_SS.dat};

\addplot [name path=upper,draw=none] table[x=averageTimes,y = conUpperMedianVals] {6agent_ipoptws_SS.dat};
\addplot [name path=lower,draw=none] table[x=averageTimes,y = conLowerMedianVals] {6agent_ipoptws_SS.dat};
\addplot [fill=ipoptws_SS_color!20] fill between[of=upper and lower];

\definecolor{ipoptran_SS_color}{rgb}{0.9290,0.6940,0.1250}
\addplot [mark=none, color=ipoptran_SS_color,line width=1pt] table[x=averageTimes,y=conMedianVals] {6agent_ipoptran_SS.dat};
\addplot [name path=upper,draw=none] table[x=averageTimes,y = conUpperMedianVals] {6agent_ipoptran_SS.dat};
\addplot [name path=lower,draw=none] table[x=averageTimes,y = conLowerMedianVals] {6agent_ipoptran_SS.dat};
\addplot [fill=ipoptran_SS_color!20] fill between[of=upper and lower];

\definecolor{proxws_color}{rgb}{1,0,0}
\addplot [mark=none, color=proxws_color,line width=1pt] table[x=averageTimes,y=conMedianVals] {6agent_proxws.dat};

\addplot [name path=upper,draw=none] table[x=averageTimes,y = conUpperMedianVals] {6agent_proxws.dat};
\addplot [name path=lower,draw=none] table[x=averageTimes,y = conLowerMedianVals] {6agent_proxws.dat};
\addplot [fill=proxws_color!20] fill between[of=upper and lower];

\draw[dashed, very thick, color=black] (axis cs:5.73827373350000, 1e-15) -- (axis cs:5.73827373350000, 1e13);

\end{loglogaxis}
\end{tikzpicture}}
\caption{Constraint violation (6 agents).}
\end{subfigure}
\caption{Convergence of the sequential convex programming approach proposed in Sec.~\ref{sec: proposed SCP} and benchmark methods with warm-starting and random initialization for 100 Monte Carlo simulations. The objective function corresponds to the objective function in optimization~\eqref{opt: min_time_relaxed_v1} after integrating---from the initial state \(\overline{\bm{x}}_0\)---the dynamics in optimization~\eqref{opt: min_time_relaxed_v1} using the computed input sequences with \texttt{ode45}, and the constraint violation corresponds to~\eqref{constraint_violation}. In the figure, the black dashed vertical line represents the median time spent obtaining the warm-starting using Algorithm~\ref{alg: filter}. The solid lines represent the median values of the simulations, while the shaded areas indicate the interquartile range, with the lower bound at the first quartile (0.25 quantile) and the upper bound at the third quartile (0.75 quantile). 
}

\label{fig: Convergence_Obj_Con}
\end{figure}

\subsection{Simulation Results}

We illustrate the convergence performance of the proposed SCP approach and the benchmark methods in Fig.~\ref{fig: Convergence_Obj_Con}. Note that, instead of using both the computed states and the input sequences, we perform post‐processing by integrating---from the initial state \(\overline{\bm{x}}_0\)---the dynamics in optimization~\eqref{opt: min_time_relaxed_v1} using only the computed input sequences with \texttt{ode45}. We denote the computed input sequence by $\{\tilde{\bm{\eta}}_k\}_{k=1}^{N-1}$ and the state sequence obtained by integrating from $\overline{\bm{x}}_0$ under these inputs by $\{\tilde{\bm{\xi}}_k\}_{k=2}^{N}$. Because this integration enforces the dynamics constraints, the violations reported in Fig.~\ref{fig: Convergence_Obj_Con} reflect only input constraint violations and the terminal state violation at $\tilde{\bm{\xi}}_N$, \ie,
\begin{equation}
\label{constraint_violation}
\begin{bmatrix}
    0_{n_x}^\top & 1 & 0 
\end{bmatrix} \tilde{\bm{\xi}}_N
+ \left\lVert
\begin{bmatrix}
    I_{n_x} & 0_{n_x \times 2} 
\end{bmatrix}
\left( \tilde{\bm{\xi}}_N - \overline{\bm{x}}_f \right)
\right\rVert_1 + \sum_{k=1}^{N-1} \left\{ \norm{\max\{\tilde{\bm{\eta}}_k - \overline{\bm{u}}_{\max},0_{n_u+1} \}}_1 + \norm{\max\{ -\tilde{\bm{\eta}}_k + \overline{\bm{u}}_{\min} , 0_{n_u+1} \}}_1 \right\}.
\end{equation}
In the following convergence comparison, we compare the methods at the point where the objective functions converge to the same level and the median constraint violations fall below \(10^{-2}\). Under random initialization, the proposed SCP approach reduces computation time by up to a factor of 20 to 30 compared with \texttt{IPOPT} using either multiple-shooting or single-shooting. The results also demonstrate that our proposed SCP approach is scalable from the two-agent to the six-agent case: each objective function curve converges within \(15\,\mathrm{s}\), whereas the curves for \texttt{IPOPT}---under both shooting schemes---have not converged or are only beginning to converge by the termination time. However, the proposed SCP approach gets stuck in suboptimal solutions, and \texttt{IPOPT} with multiple-shooting also gets stuck in the two-agent and six-agent cases until the termination time.

The results in Fig.~\ref{fig: Convergence_Obj_Con} also show that the warm‑starting strategy proposed in Section~\ref{sec: proposed warm-starting} improves the convergence performance of all methods. The proposed SCP approach with warm‑starting achieves up to a \(48\%\) lower objective function value than its randomly initialized version. Compared with \texttt{IPOPT}, the proposed SCP approach combined with the warm-starting strategy reduces computation time by up to two orders of magnitude. Furthermore, compared with random initialization, the warm-starting strategy consistently helps all methods reach a lower objective value and a faster reduction of constraint violations. Thus, our warm-starting strategy is shown to be compatible not only with the proposed SCP approach, but also with the benchmark methods.



\begin{figure}[!htp]

\begin{tikzpicture}[baseline]
\pgfplotslegendfromname{sharedlegend_2}
\end{tikzpicture}
\centering

\begin{subfigure}{0.47\columnwidth}
    \centering
    \pgfplotsset{every tick label/.append style={font=\scriptsize}}
\def\ratio{0.85\linewidth}
\begin{tikzpicture}
  \begin{axis}[
    ymin=0, ymax=12,
    xmin=0, xmax=12,
    width=\ratio,
    height=\ratio,
    xlabel=\footnotesize Time (s),
    ylabel=\footnotesize Distance to Obstacle 1 Center (\(\mathrm{m}\)),
    xtick={0,6,12},
    ytick={0,2, 6, 12},
    legend style={
    at={(rel axis cs:1.02,0.5)},  
    anchor=west,                  
    legend columns=2,             
    /tikz/column 2/.style={column sep=10pt}, 
    row sep=4pt,                  
    draw=black,
    font=\scriptsize
  },
    legend image post style={line width=1.5pt},
  ]

    \addplot[solid, thick, color = color1]
      table[x=time, y=d1] {distance_o1.dat};

    \addplot[
      only marks,
      mark=o,
      mark options={draw= color1, fill=none, scale=1}
    ]
    coordinates {
    (0, 6.70820393249937)
(1.00338571591280, 6.54209299375965)
(2.12549155801555, 5.29579514354369)
(3.34333191796246, 2.803663488049345)
(4.92891086502592, 3.17039765204507)
(7.30726195334066, 7.32019828174113)
(9.25584877321507, 10.02886863004357)
(11.8491363512125, 10.81665333097342)
    };

    \addplot[solid, thick, color = color2]
      table[x=time, y=d2] {distance_o1.dat};

    \addplot[
      only marks,
      mark=o,
      mark options={draw= color2, fill=none, scale=1}
    ]
    coordinates {
(0, 10.81665382639197)
(1.00338571591280, 10.67100442088402)
(2.12549155801555, 9.53171544344895)
(3.34333191796246, 6.89565238274833)
(4.92891086502592, 3.68811679068101)
(7.30726195334066, 3.06538192437980)
(9.25584877321507, 5.72278809950423)
(11.8491363512125, 6.70820343138422)
    };

    \addplot[solid, thick, color = color3]
      table[x=time, y=d3] {distance_o1.dat};

    \addplot[
      only marks,
      mark=o,
      mark options={draw= color3, fill=none, scale=1}
    ]
    coordinates {
(0, 10.81665382639197)
(1.00338571591280, 10.72296534616409)
(2.12549155801555, 9.93013909896303)
(3.34333191796246, 7.81257764910183)
(4.92891086502592, 4.73694994879311)
(7.30726195334066, 2.83645252896548)
(9.25584877321507, 5.62324055430728)
(11.8491363512125, 6.70820326196396)
    };

    \addplot[solid, thick, color = color4]
      table[x=time, y=d4] {distance_o1.dat};

    \addplot[
      only marks,
      mark=o,
      mark options={draw= color4, fill=none, scale=1}
    ]
    coordinates {
    (0, 10.81665382639197)
(1.00338571591280, 10.71836269918872)
(2.12549155801555, 10.01111627629861)
(3.34333191796246, 8.47221323005059)
(4.92891086502592, 5.90504875406018)
(7.30726195334066, 4.34509031533792)
(9.25584877321507, 6.11332564406093)
(11.8491363512125, 6.70820373973571)
};

    \addplot[solid, thick, color = color5]
      table[x=time, y=d5] {distance_o1.dat};

    \addplot[
      only marks,
      mark=o,
      mark options={draw= color5, fill=none, scale=1}
    ]
    coordinates {
(0, 9)
(1.00338571591280, 8.98093197738329)
(2.12549155801555, 8.79342800954392)
(3.34333191796246, 8.04163457055572)
(4.92891086502592, 5.89250421000320)
(7.30726195334066, 3.10065811781886)
(9.25584877321507, 2.32097942956816)
(11.8491363512125, 2.99999968787353)
  };

    \addplot[solid, thick, color = color6]
      table[x=time, y=d6] {distance_o1.dat};

    \addplot[
      only marks,
      mark=o,
      mark options={draw= color6, fill=none, scale=1}
    ]
    coordinates {
    (0, 3)
(1.00338571591280, 2.977727297495882)
(2.12549155801555, 2.698325225806842)
(3.34333191796246, 2.117568173803255)
(4.92891086502592, 2.914076402246926)
(7.30726195334066, 6.13134843110233)
(9.25584877321507, 8.35414049155855)
(11.8491363512125, 8.99999971599629)
};

\draw[dashed, thick, color=black] (axis cs:0, 2) -- (axis cs:12, 2);

\end{axis}
\end{tikzpicture}
    \caption{Distance between each agent and the center of obstacle~1 over time.}
\end{subfigure}
\hfill
\begin{subfigure}{0.47\columnwidth}
    \centering
    \pgfplotsset{every tick label/.append style={font=\scriptsize}}
\def\ratio{0.85\linewidth}
\begin{tikzpicture}
  \begin{axis}[
    ymin=0, ymax=12,
    xmin=0, xmax=12,
    width=\ratio,
    height=\ratio,
    xlabel=\footnotesize Time (s),
    ylabel=\footnotesize Distance to Obstacle 2 Center (\(\mathrm{m}\)),
    xtick={0,6,12},
    ytick={0,1.5,6,12},
    legend style={
    at={(rel axis cs:1.02,0.5)},  
    anchor=west,                  
    legend columns=2,             
    /tikz/column 2/.style={column sep=10pt}, 
    row sep=4pt,                  
    draw=black,
    font=\scriptsize
  },
    legend image post style={line width=1.5pt},
  ]


    \addplot[solid, thick, color = color1]
      table[x=time, y=d1] {distance_o2.dat};

    \addplot[
      only marks,
      mark=o,
      mark options={draw= color1, fill=none, scale=1}
    ]
    coordinates {
(0, 7.61577310586391)
(1.00338571591280, 7.44412929880903)
(2.12549155801555, 6.16476501687231)
(3.34333191796246, 3.61499478008273)
(4.92891086502592, 3.44240262566823)
(7.30726195334066, 7.22435264264797)
(9.25584877321507, 9.59729637433777)
(11.8491363512125, 10.29562969919216)
    };

    \addplot[solid, thick, color = color2]
      table[x=time, y=d2] {distance_o2.dat};

    \addplot[
      only marks,
      mark=o,
      mark options={draw= color2, fill=none, scale=1}
    ]
    coordinates {
(0, 5.83095189484530)
(1.00338571591280, 5.68802915275076)
(2.12549155801555, 4.58262958927925)
(3.34333191796246, 2.26343855140627)
(4.92891086502592, 2.45316896037828)
(7.30726195334066, 6.82856521850570)
(9.25584877321507, 10.3166371100763)
(11.8491363512125, 11.4017536931851)
    };

    \addplot[solid, thick, color = color3]
      table[x=time, y=d3] {distance_o2.dat};

    \addplot[
      only marks,
      mark=o,
      mark options={draw= color3, fill=none, scale=1}
    ]
    coordinates {
(0, 10.2956301409870)
(1.00338571591280, 10.1948069962139)
(2.12549155801555, 9.35274664620497)
(3.34333191796246, 7.10651257384719)
(4.92891086502592, 3.66861061745294)
(7.30726195334066, 3.13675057619573)
(9.25584877321507, 6.48943043122074)
(11.8491363512125, 7.61577248734234)
    };

    \addplot[solid, thick, color = color4]
      table[x=time, y=d4] {distance_o2.dat};

    \addplot[
      only marks,
      mark=o,
      mark options={draw= color4, fill=none, scale=1}
    ]
    coordinates {
(0, 5.83095189484530)
(1.00338571591280, 5.73567422136994)
(2.12549155801555, 5.06652694999968)
(3.34333191796246, 3.81896052710021)
(4.92891086502592, 3.48832201321240)
(7.30726195334066, 7.93570471768962)
(9.25584877321507, 10.6855464575635)
(11.8491363512125, 11.4017540367561)
    };

    \addplot[solid, thick, color = color5]
      table[x=time, y=d5] {distance_o2.dat};

    \addplot[
      only marks,
      mark=o,
      mark options={draw= color5, fill=none, scale=1}
    ]
    coordinates {
(0, 5.83095189484530)
(1.00338571591280, 5.81975320411851)
(2.12549155801555, 5.71606684513722)
(3.34333191796246, 5.35852800341148)
(4.92891086502592, 5.00873311185019)
(7.30726195334066, 6.75099118100725)
(9.25584877321507, 7.26804966735192)
(11.8491363512125, 7.61577288742888)
    };

    \addplot[solid, thick, color = color6]
      table[x=time, y=d6] {distance_o2.dat};

    \addplot[
      only marks,
      mark=o,
      mark options={draw= color6, fill=none, scale=1}
    ]
    coordinates { 
    (0, 7.61577310586391)
(1.00338571591280, 7.63406236034598)
(2.12549155801555, 7.59630183789685)
(3.34333191796246, 6.92311674083735)
(4.92891086502592, 5.23415754613267)
(7.30726195334066, 4.30808807603345)
(9.25584877321507, 5.37703431902207)
(11.8491363512125, 5.83095163243266)
    };

\draw[dashed, thick, color=black] (axis cs:0, 1.5) -- (axis cs:12, 1.5);

    \end{axis}
    \end{tikzpicture}
    \caption{Distance between each agent and the center of obstacle~2 over time.}
\end{subfigure}

\bigskip

\begin{subfigure}[t]{0.47\columnwidth}
    \centering
    \pgfplotsset{every tick label/.append style={font=\scriptsize}}

\def\ratio{0.85\linewidth}

\begin{tikzpicture}
  \begin{axis}[
    ymin=0, ymax=3.2,
    xmin=0, xmax=12,
    width=\ratio,
    height=\ratio,
    xlabel=\footnotesize Time (s),
    ylabel=\footnotesize Velocity Magnitude (\(\mathrm{m/s}\)),
    xtick={0,6,12},
    ytick={0,1,2,3},
   legend style={
    at={(rel axis cs:1.02,0.5)},  
    anchor=west,                  
    legend columns=2,             
    /tikz/column 2/.style={column sep=10pt}, 
    row sep=4pt,                  
    draw=black,
    font=\scriptsize
  },
    legend image post style={line width=1.5pt},
  ]
  

    \addplot[solid, thick, red]
      table[x=time, y=y1] {data_velocity_ws.dat};

    \addplot[
      only marks,
      mark=o,
      mark options={draw=red, fill=none, scale=1}
    ]
      coordinates {
        (0,0)
        (1.0034,0.5644)
        (2.1255,2.0352)
        (3.3433,2.9209)
        (4.9289,2.8021)
        (7.3073,2.6972)
        (9.2558,1.5101)
        (11.8491,5.5619e-07)
      };

    \addplot[solid, thick, color = color2]
      table[x=time, y=y2] {data_velocity_ws.dat};

    \addplot[
      only marks,
      mark=o,
      mark options={draw= color2, fill=none, scale=1}
    ]
    coordinates {
      (0,0)
      (1.00338571591280,0.535075493180443)
      (2.12549155801555,1.97788184097584)
      (3.34333191796246,2.71809218162004)
      (4.92891086502592,2.51082749540959)
      (7.30726195334066,2.71498297572843)
      (9.25584877321507,1.62219880460689)
      (11.8491363512125,8.01981776992960e-07)
    };
    \addplot[solid, thick, color = color3]
      table[x=time, y=y3] {data_velocity_ws.dat};

    \addplot[
      only marks,
      mark=o,
      mark options={draw= color3, fill=none, scale=1}
    ]
    coordinates {
      (0,0)
      (1.00338571591280,0.353820439722213)
      (2.12549155801555,1.48665906187167)
      (3.34333191796246,2.55948730607661)
      (4.92891086502592,2.89803633571843)
      (7.30726195334066,2.82224823729969)
      (9.25584877321507,1.58407600244501)
      (11.8491363512125,1.00310730182802e-06)
    };

    \addplot[solid, thick, color = color4]
      table[x=time, y=y4] {data_velocity_ws.dat};

    \addplot[
      only marks,
      mark=o,
      mark options={draw= color4, fill=none, scale=1}
    ]
    coordinates {
      (0,0)
      (1.00338571591280,0.381738704593372)
      (2.12549155801555,1.51253594059612)
      (3.34333191796246,2.69156267313300)
      (4.92891086502592,2.90079181562353)
      (7.30726195334066,2.82227100853831)
      (9.25584877321507,2.03546738958931)
      (11.8491363512125,5.98229905732512e-07)
    };

    \addplot[solid, thick, color = color5]
      table[x=time, y=y5] {data_velocity_ws.dat};

    \addplot[
      only marks,
      mark=o,
      mark options={draw= color5, fill=none, scale=1}
    ]
    coordinates {
      (0,0)
      (1.00338571591280,0.0743207051149820)
      (2.12549155801555,0.675591559772333)
      (3.34333191796246,2.00407005713121)
      (4.92891086502592,2.77180744795301)
      (7.30726195334066,2.50155971081784)
      (9.25584877321507,1.77791828071132)
      (11.8491363512125,1.05110700095227e-06)
    };

    \addplot[solid, thick, color = color6]
      table[x=time, y=y6] {data_velocity_ws.dat};

    \addplot[
      only marks,
      mark=o,
      mark options={draw= color6, fill=none, scale=1}
    ]
    coordinates {
      (0,0)
      (1.00338571591280,0.407688120632277)
      (2.12549155801555,1.68846159929542)
      (3.34333191796246,2.66593556852059)
      (4.92891086502592,2.32736832170030)
      (7.30726195334066,2.00380852019959)
      (9.25584877321507,1.49088845348886)
      (11.8491363512125,6.77369121462759e-07)
    };

\draw[dashed, thick, color=black] (axis cs:0, 3) -- (axis cs:12, 3);

  \end{axis}
\end{tikzpicture}
    \caption{Velocity Magnitudes of six agents over time.}
\end{subfigure}
\hfill
\begin{subfigure}[t]{0.47\columnwidth}
    \centering
    \pgfplotsset{every tick label/.append style={font=\scriptsize}}
\def\ratio{0.85\linewidth}
\begin{tikzpicture}
  \begin{axis}[
    ymin=0, ymax=18,
    xmin=0, xmax=12,
    width=\ratio,
    height=\ratio,
    xlabel=\footnotesize Time (s),
    ylabel=\footnotesize Inter-Agent Distance (\(\mathrm{m}\)),
    xtick={0,6,12},
    ytick={0,1,6,12,18},
   legend style={
    at={(rel axis cs:1.05,0.5)},  
    anchor=west,                  
    legend columns=1,             
    /tikz/column 2/.style={column sep=8pt}, 
    row sep=2pt,                  
    draw=black,
    font=\scriptsize
  },
    legend image post style={line width=1.5pt},
  ]

\addplot[solid, thin, color = black]
    table[x=time, y = d_min] {inter_agent_v.dat};


\addplot[solid, thin, color = black]
    table[x=time, y = d_max] {inter_agent_v.dat};


\addplot [name path=upper,draw=none] table[ x = time, y = d_max ] {inter_agent_v.dat};
\addplot [name path=lower,draw=none] table[ x = time, y = d_min ] {inter_agent_v.dat};
\addplot [fill=black!20] fill between[of=upper and lower];

\draw[dashed, thick, color=black] (axis cs:0, 1) -- (axis cs:12, 1);
    
  \end{axis}
\end{tikzpicture}

  
    \caption{Range of inter-agent distances over time across all agent pairs.}
\end{subfigure}

\caption{Computed trajectories for six agents using the sequential convex programming approach proposed in Sec.~\ref{sec: proposed SCP} with warm-starting, corresponding to the realization with the median objective function value among 100 Monte Carlo simulations terminated at 15 seconds. The trajectories are computed by integrating---from the initial state \(\overline{\bm{x}}_0\)---the dynamics in optimization~\eqref{opt: min_time_relaxed_v1} using the computed input sequences with \texttt{ode45}. The horizontal dashed line in each figure corresponds to the respective constraint.
}

\label{fig: SIX_WS}
\end{figure}
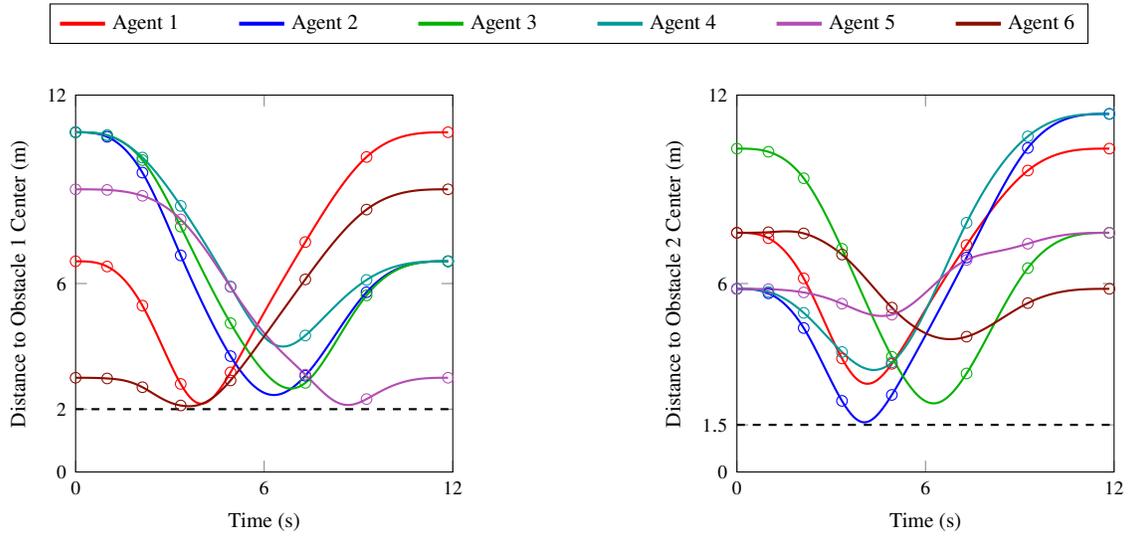
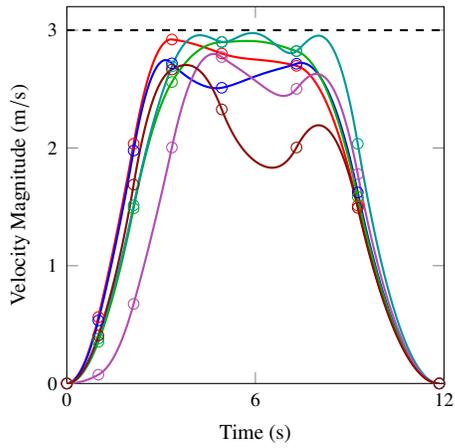
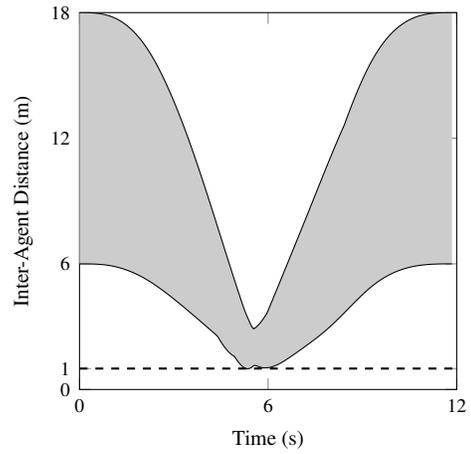

\begin{figure}
  \centering

  \begin{subfigure}{0.47\columnwidth}
    \centering
    
  \newlength{\picw}\setlength{\picw}{0.85\linewidth}  
  \newlength{\pich}\setlength{\pich}{0.6\picw}    

  \begin{tikzpicture}[x=1cm,y=1cm]                
    \node[inner sep=0pt, anchor=south west] (img) at (0,0)
          {\includegraphics[width=\picw]{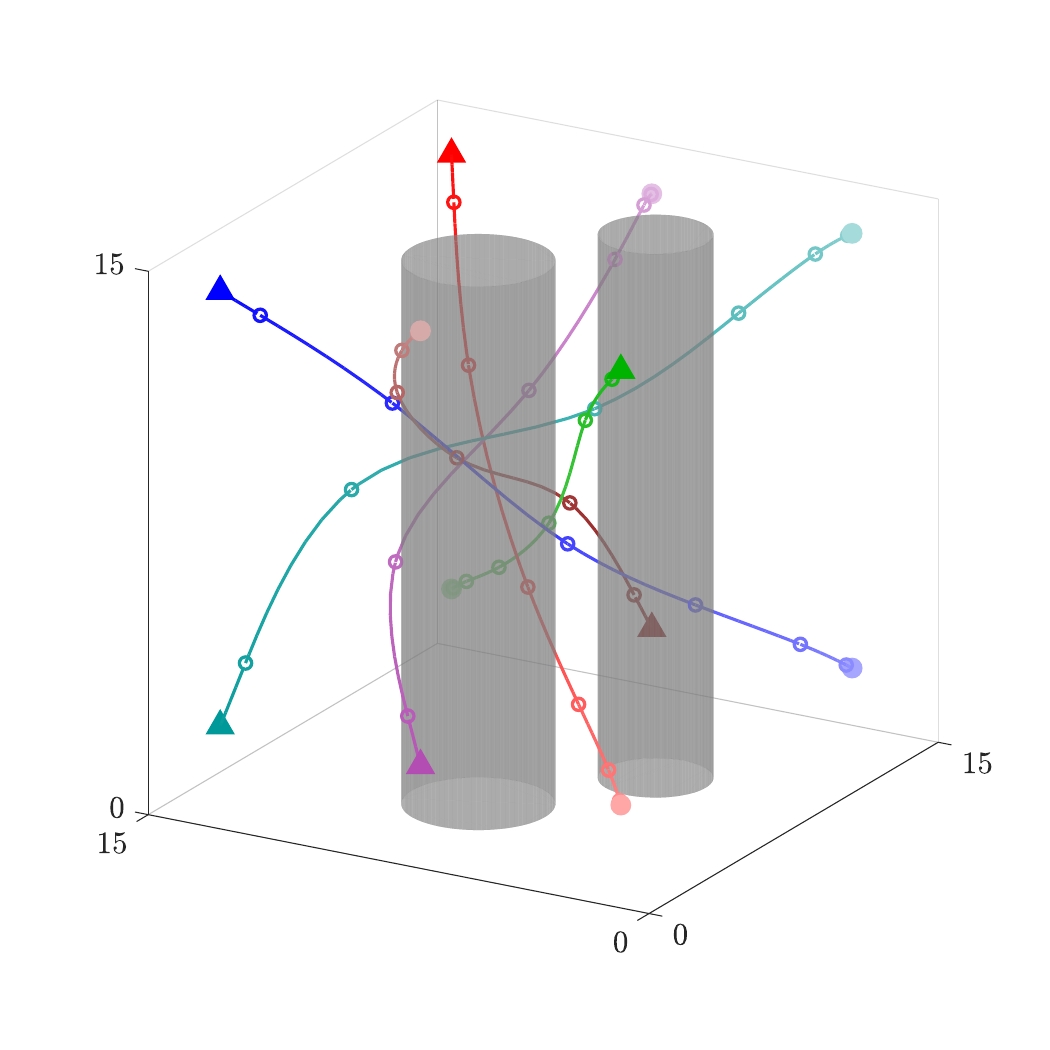}};

    \node[font=\footnotesize]
          at ({\picw/1.25}, 0.9) {$X\;(\mathrm{m})$};

    \node[font=\footnotesize, anchor=south west]
          at ({\picw/2.8}, 0.5) {$Y\;(\mathrm{m})$};

    \node[font=\footnotesize, rotate=90]
          at (0.5, {\pich/1.2}) {$Z\;(\mathrm{m})$};
          
  \end{tikzpicture}
    \caption{Warm-Starting via Algorithm~\ref{alg: filter}}
    \label{plot: trajectory with warm-starting}
    \end{subfigure}
    \hfill
  \begin{subfigure}{0.48\columnwidth}
  \centering
  \newlength{\picwid}
  \setlength{\picwid}{0.85\linewidth}
  \newlength{\picheight}
  \setlength{\picheight}{0.6\picwid}         

  \begin{tikzpicture}[x=1cm,y=1cm]
    \node[inner sep=0pt, anchor=south west] (img) at (0,0)
      {\includegraphics[width=\picwid]{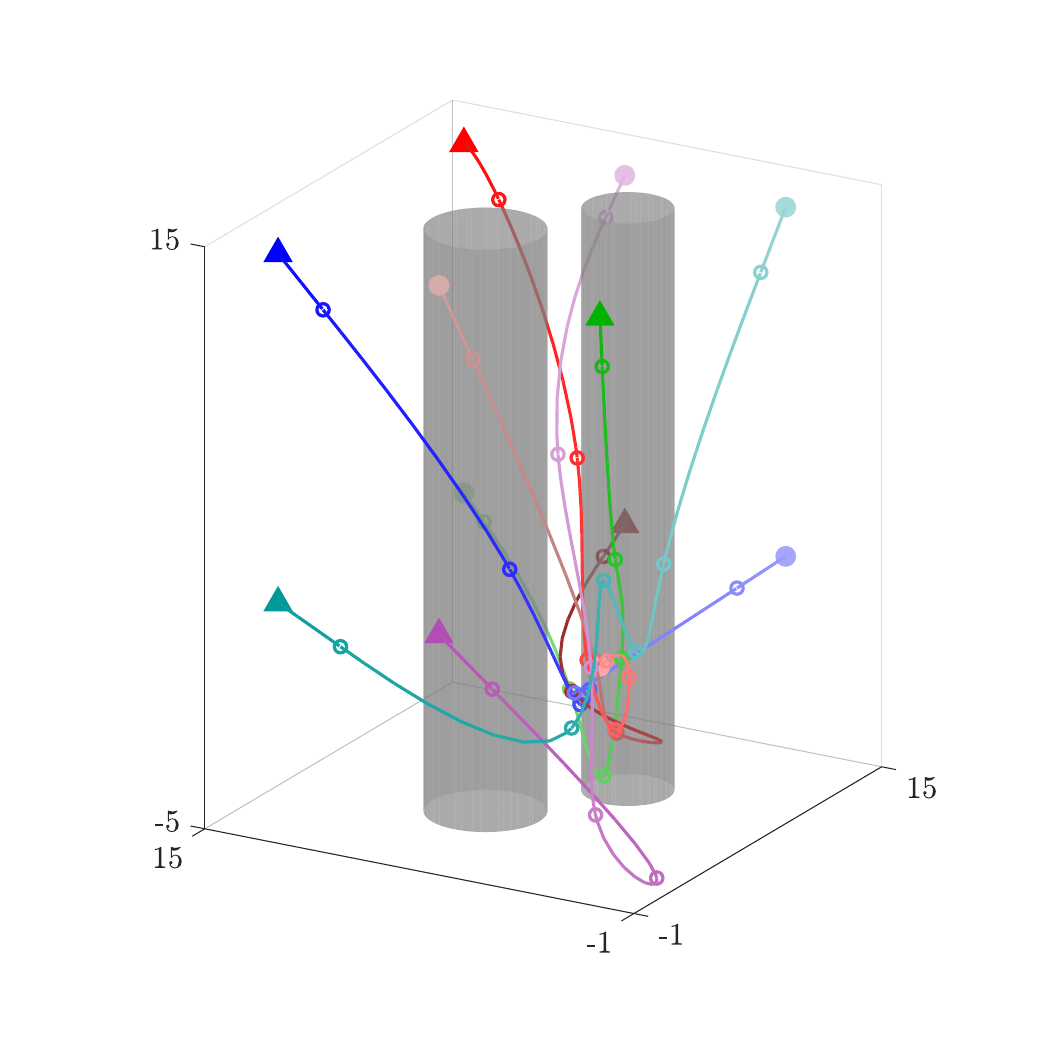}};

    \node[font=\footnotesize]    at ({\picwid/1.25}, 0.9)
      {$X\;(\mathrm{m})$};
    \node[font=\footnotesize, anchor=south west]
      at ({\picwid/2.8}, 0.5) {$Y\;(\mathrm{m})$};
    \node[font=\footnotesize, rotate=90]  at (0.8, {\picheight/1.2})
      {$Z\;(\mathrm{m})$};
  \end{tikzpicture}
  \caption{Random initialization}
  \label{plot: trajectory with random initialization}
  \end{subfigure}

\caption{Computed trajectories for six agents generated by the sequential convex programming approach proposed in Sec.~\ref{sec: proposed SCP} with warm-starting and random initialization, respectively. The two plots correspond to realizations with the median objective function value among 100 Monte Carlo simulations terminated at 15 seconds. The trajectories are computed by integrating---from the initial state \(\overline{\bm{x}}_0\)---the dynamics in optimization~\eqref{opt: min_time_relaxed_v1} using the computed input sequences with \texttt{ode45}. All agents start from the circle (initial position) and move to the triangle (final position). Along each trajectory, the discretization points include the unfilled circles, as well as the initial circle and final triangle. As each trajectory approaches its end, its color becomes progressively darker.}

\label{fig: SIX_wr}
    
\end{figure}

We illustrate the computed trajectories---including velocity magnitudes, distances to obstacles, inter-agent distances, and position trajectories---in Fig.~\ref{fig: SIX_WS} and Fig.~\ref{fig: SIX_wr}. These plots correspond to the trajectories that yield the median objective function value among the 100 Monte Carlo simulations terminated at $15\,\mathrm{s}$. As shown in Fig.~\ref{fig: SIX_WS}, our proposed SCP approach ensures continuous-time constraint satisfaction almost everywhere---that is, not only at the discretization points, but also between them. As shown in Fig.~\ref{plot: trajectory with warm-starting}, the trajectories obtained with warm-starting are smoother than those in Fig.~\ref{plot: trajectory with random initialization}. The latter violate the box-constraints at the termination time.

\section{Conclusion}

We proposed a framework for optimizing multiagent quadrotor trajectories that combines (i) sequential convex programming with system dynamics that embed continuous-time constraints, and (ii) filtering-based warm-starting using constraint-aware particle filtering. This framework ensures continuous-time constraint satisfaction. Compared with random initialization, the warm-starting strategy consistently helps both the proposed SCP approach and benchmark methods reach a lower objective value and a faster reduction in constraint violations.



However, the current work still has limitations. For example, it remains unclear how to efficiently adjust the parameters in Algorithm~\ref{alg: filter} for different problems. In addition, it is also unclear how the current can impact the performance of multiagent model predictive control. In future work, we aim to address these limitations, as well as explore applications beyond quadrotor systems.

\bibliography{sample}

\end{document}